\documentclass[12pt,russian,english]{article} 
\usepackage[T2A]{fontenc}  
\usepackage{babel,amsfonts,latexsym,amssymb,amsmath}

\def \al{\alpha}
\def \be{\beta}

\def \dl{\delta}
\def \ep{\varepsilon}

\def \ka{\varkappa}

\def \ph{\varphi}
\def \om{\omega}

\def \Om{\Omega}

\def \operatorname#1{\mathop{\rm #1}}

\def\div{\operatorname{div}}
\def\osc{\operatorname{osc}}
\def\osc2{\operatorname{osc^2}}
\def\dist{\operatorname{dist}}

\def\cd{\partial}

\def\Q0{Q(x_0,t_0,R)}
\def\0{{x_0,t_0,R}}
\def\build#1_#2{\mathrel{\mathop{\kern 0pt#1}\limits_{#2}}}

\newcommand{\pint}{\mathop{\int\limits{\hspace{-4mm}}-}\limits}
\newtheorem{theorem}{Theorem}[section]
\newtheorem{definition}{Definition}[section]
\newtheorem{counterexample}{Example}[section]

\newtheorem{lemma}{Lemma}[section]

\title{The Local Regularity Theory \\ for the Navier--Stokes Equations \\ Near the Boundary}
\author{G.A.~Seregin and T.N.~Shilkin\thanks{This work is  supported by RFBR grant  11-01-00324}}

\begin{document}

\maketitle

\abstract{This is an expository paper on 
the theory of local regularity for weak solutions to the non-stationary 3D Navier-Stokes equations near the boundary of a domain.}

\section{Introduction}

\bigskip
\medskip
The main problem of the modern mathematical hydrodynamics is the global well-posedness of the 3D Navier-Stokes equations, i.e. the global existence of a unique solution, corresponding to a given smooth divergent-free initial data.
There are two main directions 
in the study of this problem. One can try to improve local well-posedness results, see classical papers \cite{Leray} and \cite{Kiselev_Ladyzhenskaya} and many others on local well-posedness, or to show that a global week solution, introduced essentially in \cite{Leray} and  \cite{Hopf} and called the weak Leray-Hopf solution,  to the corresponding initial boundary value problem is in fact unique. On the other hand, the second aim can be achieved by proving regularity of weak solutions. Indeed, it is well-known, since the celebrated paper \cite{Leray} has been published, that smoothness of weak solutions implies their uniqueness in the class of weak Leray-Hopf solutions. Here, we are going to discuss regularity of weak solutions keeping in mind that it is one of possible ways to attack the main problem on the global well-posedness. Our approach is quite  typical for PDE's theory and in a sense local. The latter means that we have a solution to the Navier-Stokes system with a ``finite energy'' in a canonical parabolic domain and try to show that it is smoother in subdomains. The result depends on assumptions imposed on the pressure. Our choice of the class for the pressure is motivated by the linear theory and reflects
the fact the whole Navier-Stokes (or Stokes) problem is not quite local because of the incompressibility condition.

Given a space-time point $z_0=(x_0,t_0) $, the canonical domain is going to be  
  a parabolic cylinder $Q(z_0, R):= B(x_0,R)\times ]t_0-R^2, t_0[$ if the local interior regularity is studied
  , or parabolic half-cylinder $Q^+(z_0,R):=B^+(x_0,R)\times ]t_0-R^2,t_0[$ if the local boundary regularity is under consideration. Here $B(x_0,R)$ denotes a ball in $\Bbb R^3$ of  radius $R$ centered at a point $x_0$, and $B^+(x_0,R)$ is a half-ball $B(x_0,R)\cap (x_0+ \Bbb R^3_+)$, and $\Bbb R^3_+:=\{ ~x=(x_1,x_2,x_3)\in \Bbb R^3~|~ x_3>0~\}$.

Since the Navier-Stokes system is invariant under the scaling  transformation
\begin{equation}
u^R(y,s) = Ru(x_0+Ry, t_0+R^2s), \qquad p^R(y,s) = R^2 p(x_0+Ry, t_0+R^2s),
 \label{Scale_Invariance}
\end{equation}
the problem of local regularity  of weak solutions in a neighborhood of a point $z_0$ can be  reduced to a model problem in some fixed  domain, say, $Q=B\times ]-1,0[$ in the internal case or $Q^+= B^+\times ]-1,0[$ in the boundary case. Here, $B$ is the unite ball of $\Bbb R^3$ centered at the origin and $B^+:=B\cap \Bbb R^3_+$.

In contrast to equations of parabolic type,  the local smoothing for the Navier-Stokes  system has some special features. In particular, in the local setting, 
weak solutions  may be not infinitely smooth in  subdomains despite  the right-hand side is infinitely smooth  there. So, there might be a limiting smoothness which can be achieved locally.  The roots of this phenomena lay in the linear theory which we discuss detailly in Section \ref{Linear_Theory}.

Typical results in the known local regularity theory have the form of the 
so-called $\ep$-regularity conditions. For the Navier-Stokes system in the canonical domains, $\ep$-regularity  conditions ensure  H\" older continuity of the velocity field  around the  origin  provided  a certain integral quantity of a solution over the above domain is sufficiently  small.

The first results on the local regularity for the 3D Navier-Stokes equations belong to Scheffer \cite{Scheffer}. Scheffer considered a special class of weak solutions to the Cauchy problem  that satisfy  a {\it local energy inequality}. Motivated by this observation, later on, Caffarelly-Kohn-Nirenberg introduced the so-called {\it suitable weak solutions} that are just solutions to the Navier-Stokes system with certain reasonable properties. That was a great step towards a complete local setting. By a definition, a suitable weak solution $u$ and $p$ is such that the velocity $u$ belongs to the energy class, pressure is an integrable function (with a certain exponent of integrability determined by  the linear theory and by the integrability of the convective term),   $u$ and $p$ are assumed to satisfy the Navier-Stokes system in the sense of distributions and  the local energy inequality. We will explore the definition of suitable weak solutions introduced in \cite{Lin}, see also \cite{Ladyzhenskaya_Seregin}:

\begin{definition}\label{Definition_1}
We say that a pair of functions $u$ and $p$ is   {\it a suitable weak solution} to the Navier-Stokes system in $Q$ if
\begin{itemize}
\item $u\in L_{2, \infty}(Q)\cap W^{1,0}_2(Q)$, $p\in L_{\frac 32}(Q)$
\item $u$ and $p$ satisfy the Navier-Stokes system in $Q$ in the sense of distributions
\item for a.a. $t\in ]-1,0[$, the pair $u$ and $p$ satisfies the local energy inequality in $Q$
$$
\gathered
\int\limits_{B} \zeta(x,t) |u(x,t)|^2~dx \ + \ 2\int\limits^t_{-1}\int\limits_{B} \zeta |\nabla u|^2~dxdt \ \le \\ \le \
\int\limits^t_{-1}\int\limits_{B}  |u|^2~ \left(\cd_t \zeta +\Delta \zeta \right)~dxdt \ + \
\int\limits^t_{-1}\int\limits_{B}  u\cdot\nabla \zeta ~\left(|u|^2 + 2p \right)~dxdt
\endgathered
$$
for any non-negative test function $\zeta\in C^\infty(\mathbb R^3\times\mathbb R)$ 
vanishing  near the parabolic boundary  $\cd'Q:= (\cd B\times ]-1,0[)\cup (B\times \{ t=-1\})$.
\end{itemize}
\end{definition}
Here we denote by $L_s(Q)$ the Lesbegue space of functions integrable over $Q$ with the exponent $s\in [1, +\infty]$; \ $W^{1,0}_s(Q):=\{ u\in L_s(Q)~|~ \nabla u\in L_s(Q)\}$, where $\nabla u$ denotes  the  gradient of $u$ with respect to spatial variables; $L_{2,\infty}(Q):= L_\infty(-1, 0; L_2(B))$.

This class of local solutions appears in the global setting. Indeed, it is not so difficult to show that, in a given space-time domain, among all weak Leray-Hopf solutions, corresponding to a given initial data, there exists at least one that is a suitable weak solution in any parabolic cylinder $Q(z_0,R)$ belonging to the space-time domain, see \cite{CKN} and \cite{Ladyzhenskaya_Seregin}.

The significant contribution into the local regularity theory for the Navier-Stokes equation has been made by Caffarelli-Kohn-Nirenberg.  They showed that the set of all singular points is very small in the following sense: the one dimensional parabolic Hausdorff  measure of this set is equal to zero. This is a consequence of the Caffarelli-Kohn-Nirenberg $\varepsilon$-regularity condition reading that 
 there exists an absolute constant $\ep_0$ such that, for any  suitable weak solution   $u$ and $p$ in $Q$  with
\begin{equation}
\sup\limits_{r<1} \  \frac 1r~\int\limits_{Q(r)} |\nabla u|^2~dxdt \ < \ \ep_0,
\label{CKN_condition}
\end{equation}
the velocity field $u$ is essentially bounded near the origin (actually it is H\" older continuous, as it was shown later in  \cite{Ladyzhenskaya_Seregin}). Here and in what follows we denote $Q(r):=Q(0, r)$, $Q^+(r):=Q^+(0,r)$ etc.

Among various $\ep$-regularity conditions,  see, for example, papers \cite{CKN}, \cite{Lin}, \cite{Ladyzhenskaya_Seregin},  \cite{Seregin_Congress}, and many others, we would like to point out the following one:
 there exists an absolute constant $\ep_1>0$ such that,  for any  suitable weak solution $u$ and $p$ in $Q$, satisfying in addition the inequality
\begin{equation}
\int\limits_{Q} \Big(|u|^3 + |p|^{\frac 32}\Big)~dxdt \ < \ \ep_1,
\label{Basic_Condition_Internal}
\end{equation}
the velocity field $u$ is H\"older continuous in the completion of the set $ Q(\frac 12)$. In the present paper,  condition \eqref{Basic_Condition_Internal} and its boundary analogue are called  {\it the basic $\ep$-regularity conditions}.
The basic $\ep$-regularity condition is remarkable as many other $\ep$-regularity conditions, including CKN-condition \eqref{CKN_condition}, can be derived from this one, see, for example, \cite{CKN}, \cite{Lin}, and \cite{Ladyzhenskaya_Seregin}. In what follows,  we are going to deal with the boundary analogue of condition \eqref{Basic_Condition_Internal}.


Now, let us review known results on local regularity up to the spatial boundary for weak solutions to the Navier-Stokes system.  
In this case, we complement the Navier-Stokes equation with the non-slip boundary condition for the velocity field. We focus ourselves mostly on the explanation of what happens with weak solutions to the Navier-Stokes system in the canonical domain $Q^+$ with the boundary condition on the flat part of semi-ball $B^+$, i.e.,
$$
u|_{x_3=0} = 0.
$$

An analogue of $\varepsilon$-regularity condition (\ref{CKN_condition}) for boundary points has been proven in \cite{Seregin_JMFM}. The proof is typical for PDE's system and based contradiction arguments. 
Later on, in \cite{Seregin_Aa},  the same sufficient regularity condition  was proved directly, which made it possible in principle to estimate the size of all constants in the corresponding assumptions. 

The latest version of the definition  of suitable weak solutions is as follows.

\begin{definition} \label{Definition_2}
We say that a pair of functions $u$ and $p$ is   {\it a boundary suitable weak solution} to the Navier-Stokes system in $Q^+$ if
\begin{itemize}
\item $u\in L_{2, \infty}(Q^+)\cap W^{1,0}_2(Q^+)$, $p\in L_{\frac 32}(Q^+)$
\item $u|_{x_3=0}=0$ in the sense of traces
\item $u$ and $p$ satisfy the Navier-Stokes system in $Q^+$ in the sense of distributions
\item for a.a. $t\in ]-1,0[$, the pair $u$ and $p$ satisfies the local energy inequality in $Q^+$
$$
\gathered
\int\limits_{B^+} \zeta(x,t) |u(x,t)|^2~dx \ + \ 2 \int\limits^t_{-1}\int\limits_{B^+} \zeta |\nabla u|^2~dxdt \ \le \\ \le \
\int\limits^t_{-1}\int\limits_{B^+}  |u|^2~ \left(\cd_t \zeta +\Delta \zeta \right)~dxdt \ + \
\int\limits^t_{-1}\int\limits_{B^+}  u\cdot\nabla \zeta ~\left(|u|^2 + 2p \right)~dxdt
\endgathered
$$
for any non-negative test function $\zeta\in C^\infty(\mathbb R^3\times\mathbb R)$ 
vanishing  near the parabolic boundary  $\cd'Q$.
\end{itemize}

\end{definition}
However, in \cite{Seregin_JMFM}, \cite{Seregin_Aa}, \cite{SSS}, the definition of suitable weak solutions is different. It is supposed there in addition that
the   second spatial derivatives and the first  derivative in time  of the velocity field and the gradient of the pressure must exist as integrable functions in $Q^+$:
\begin{equation}
u\in W^{2,1}_{s,l}(Q^+), \quad p\in W^{1,0}_{s,l}(Q^+) \  \mbox{ for some  } \ s, l\in ]1,+\infty[ \ \mbox{ such that } \   \frac 3s+ \frac 2l\ge 4.
\label{Extra_Conditions}
\end{equation}
Here  $L_{s,l}(Q^+)$ is the anisotropic Lebesgue space equipped with the norm
$$
\|f\|_{L_{s,l}(Q^+)}:=
\Big(\int\limits_{-1}^0\Big(\int\limits_{B^+} |f(x,t)|^s~dx\Big)^{l/s}dt\Big)^{1/l} ,
$$
and we use the following notation for the functional spaces:
$$
\gathered
W^{1,0}_{s,l}(Q^+)\equiv L_l(-1,0; W^1_s(B^+))= \{ \ u\in L_{s,l}(Q^+): ~\nabla u \in L_{s,l}(Q^+) \ \},
\\
W^{2,1}_{s,l}(Q^+) = \{ \ u\in W^{1,0}_{s,l}(Q^+): ~\nabla^2 u, \ \cd_t u \in L_{s,l}(Q^+) \ \},
\endgathered
$$
and the following notation for the norms:
$$
\gathered
\| u \|_{W^{1,0}_{s,l}(Q^+)}= \| u \|_{L_{s,l}(Q^+)}+ \|\nabla u\|_{L_{s,l}(Q^+)},  \\
\| u \|_{W^{2,1}_{s,l}(Q^+)}= \| u \|_{W^{1,0}_{s,l}(Q^+)}+ \| \nabla^2 u \|_{L_{s,l}(Q^+)}+\|\cd_t u\|_{L_{s,l}(Q^+)}.
\endgathered
$$
 A particular choice of exponents $s$, $l$ is a matter of convenience and it can be made in various ways. For example, in \cite{Seregin_JMFM} and in \cite{SSS}, it is assumed that $s=9/8$, $l=3/2$, in \cite{Seregin_Aa} $s=l=5/4$ and $s=15/14$, $l=5/3$.
 It should be mentioned that the choice of $s$, $l$ is always motivated by the following general idea: if we take the Leray-Hopf solution $u$ of the initial-boundary value problem for the Navier-Stokes equation in some domain and  interpret the convective term $(u\cdot \nabla)u$ as a right-hand side of the initial-boundary value problem for the Stokes system then extra conditions \eqref{Extra_Conditions} follows  from the coercive estimates of the linear theory and from the uniqueness theorem for the Stokes problem in the class of weak solutions. This scheme works well in the case of interior regularity, while, for the boundary regularity case, its realization encounters some difficulties. Nevertheless, as it is shown in \cite{Seregin_ZNS370}, the scheme  works up to the boundary as well and extra assumptions \eqref{Extra_Conditions} are simply superfluous. 

 As it has been already mentioned, there is an essential difference between local interior regularity and local boundary regularity even for the Stokes system. One of  the consequences of such an observation is that the further smoothing in a neighborhood of regular points might happen  differently. In a neighborhood of an  interior regular point, a suitable weak solution has all the spatial derivatives that are H\"older continuous, while, in a neighborhood of boundary regular point, the spatial gradient is not necessary to be even bounded, see \cite{Kang 2005} and \cite{SS_ZNS385}.


Besides conditions \eqref{CKN_condition}, \eqref{Basic_Condition_Internal}, $\ep$-regularity theory provides many other sufficient conditions of local regularity of suitable weak  solutions to the Navier-Stokes equations, see, for example,  papers \cite {Kang 2004}, \cite{Seregin_UMN}, \cite{Seregin_Handbook}, \cite{Seregin_Congress}, and references in them. Most of these conditions are stated in terms of so-called {\it scale invariant functionals}. Some examples of such functionals (in the case of the flat part of the boundary) are as follows:
\begin{equation}
\gathered
A(u,r) =   \sup\limits_{t\in (-r^2,0)}~\Big(~\frac 1r ~\int\limits_{B^+(r)} |u(x,t)|^2~dx~\Big)^{1/2}, \\
C(u,r) =  \Big(~\frac 1{r^2} ~ \int\limits_{Q^+(r)} |u(x,t)|^3~dxdt~\Big)^{1/3}, \qquad
E(u,r) =  \Big(~\frac 1{r} ~ \int\limits_{Q^+(r)} |\nabla u(x,t)|^2~dxdt~\Big)^{1/2}
\endgathered
\label{Functionals}
\end{equation}
 ``Scale invariance'' means that if $F(u,r)$ is one of these functionals and $u^R$ and $p^R$ are  functions obtained from $u$ and $p$ by formulas \eqref{Scale_Invariance} with $x_0=0$, $t_0=0$,  then
$$
F(u^R,1) \ = \ F(u,R), \qquad \forall ~R>0.
$$
The functionals $A(u,r)$, $C(u,r)$, $E(u,r)$ possess the following property: boundedness of one of them, 
 i.e.
$$
\min \Big\{ ~\sup\limits_{r<1} A(u,r),  ~\sup\limits_{r<1} C(u,r),  ~\sup\limits_{r<1} E(u,r)~\Big\}  \ < \ +\infty,
$$
implies  boundedness of all others:
\begin{equation}
\max \Big\{ ~\sup\limits_{r<1} A(u,r),  ~\sup\limits_{r<1} C(u,r),  ~\sup\limits_{r<1} E(u,r), ~\sup\limits_{r<1} D(p,r)~\Big\}  \ < \ +\infty.
\label{Boundedness}
\end{equation}
Here, $D(p,r)$ is a functional
$$
D(p,r)  =  \Big(~\frac 1{r^2} ~ \int\limits_{Q^+(r)} |p(x,t)|^{3/2}~dxdt~\Big)^{2/3},
$$
which is also invariant under the scaling transformation of $p$ according to  \eqref{Scale_Invariance}.
Statement \eqref{Boundedness}  has been proven in \cite{Seregin_UMN} in the internal case. Later on, in \cite{Mikhaylov_1},
it  was generalized to the boundary case.


One of the basic principles in  the $\ep$-regularity theory for the Navier-Stokes equations reads: if at
  least one of the scale invariant functionals is  small uniformly with respect to all $r\in ]0,1[$, i.e.
\begin{equation}
   \sup\limits_{r<1} F(u,r) <\ep_0,
   \label{Uniform_smallness}
\end{equation}
   then the origin is a regular point of the velocity field $u$ (i.e. $u$ is H\" older continuous near the origin).
     This statement has been rigorously proven by many authors for various types of scale invariant functionals, see, for example, references  in \cite{Kang 2004}, \cite{GKT}, \cite{Seregin_Handbook}, \cite{Seregin_Congress}. In our paper, we shall mention only those contributions that are concerned with the boundary case. In \cite{Seregin_JMFM}, \cite{Seregin_Aa}, the boundary regularity up to the flat part of the boundary has been proven if $F(u,r)$ is one of functionals in \eqref{Functionals}.

There is one more example of scale invariant functionals, which we call the La\-dy\-zhen\-ska\-ya-Prodi-Serrin-type functional (LPS-functional):
$$
M_{s,l}(u,r) := \| u\|_{L_{s, l}(Q^+(r))} =  \Big(\int\limits_{-r^2}^0 \Big( \int\limits_{B^+(r)} |u(x,t)|^s~dx\Big)^{l/s}~dt\Big)^{1/l} , \qquad s,l\in ]1, +\infty[
$$
Comparably with functionals  defined in \eqref{Functionals} LPS functionals possess two additional properties. First, it is monotone with respect to $r$ and hence if it is finite for some particular  $r_0>0$ then it is uniformly bounded for all $r\in ]0,r_0[$. And second, if $s$, $l>1$, then LPS functionals are absolutely continuous functions of a domain, i.e.
\begin{equation}
M_{s,l}(u,r)\to 0 \qquad\mbox{as}\qquad r\to 0.
\label{Absolute_continuity}
\end{equation}
 So, $\ep$-regularity theory developed above provides a simple proof of the following conditional result for a boundary suitable weak solution $u$ to the Navier-Stokes equations in $Q^+$: if $s>3$ and
\begin{equation}
\frac 3s+\frac 2l =1 \qquad \mbox{and} \qquad M_{s,l}(u,1) <+\infty,
\label{LPS_condition}
\end{equation}
then $u$ is regular near the origin. For a different approach, we refer to  paper  \cite{Kang 2004}. 

Note that formally the first inequality in \eqref{LPS_condition} allows the following combination of parameters: $s=3$, $l=+\infty$, i.e. one can ask the question about  local regularity of solutions to the Navier-Stokes system belonging to the class $L_{3, \infty}(Q^+):=L_\infty(-1,0; L_3(B^+))$ which we call
$L_{3, \infty}$--solutions. It is known (see \cite{Kato}) that the initial boundary value problem for the Navier-Stokes equation is locally well-posed on $L_3$-space. But the method that we used to  prove regularity weak solution with finite LPS functionals in the case of $s>3$, $l<+\infty$ can not be extended to $L_{3,\infty}$ case as the functional $M_{s,l}(u,r)$ with $s=3$, $l=+\infty$ in general does not possess the property \eqref{Absolute_continuity}.

The proof of  regularity of $L_{3,\infty}$-solutions to the Navier-Stokes system  requires development of the completely  new approach which is based on the backward-in-time uniqueness for the heat operator with lower order coefficients in the compliment to a ball or even in a half space. This method has been introduced  in \cite{ESS} and then developed in \cite{ESS_Archive} and in \cite{ESS Algebra and Analysis}  in order to prove the interior regularity of $L_{3,\infty}$-solutions. Later on,  the same type of results was extended to the boundary case  of $L_{3,\infty}$-solutions. It has been shown in \cite{Seregin Math Ann} that:
$$
\begin{array}c
u \in L_{3,\infty}(Q^+) \qquad \Longrightarrow \qquad u\in C^{\al, \frac \al2}(\bar Q^+(\frac 12))
\end{array}
$$
Here $C^{\al, \frac \al2}(\bar Q^+(r))$ denotes the space of functions which are H\" older continuous with the exponent $\al>0$ with respect to the usual parabolic metric.

The paper is organized as follows. In Section \ref{Linear_Theory} we discuss the local smoothness of weak  solutions to  the linear Stokes problem near the boundary. In Section \ref{Basic_Condition} we present a proof of the basic $\ep$-regularity condition for boundary suitable weak solutions to the Navier-Stokes equations near a flat part of the boundary. Finally, in Section \ref{Curved_Boundary} we give a brief  overview of known results on the local regularity theory  for the Navier-Stokes equations in a domain with a curvilinear  boundary.



\section{Linear Theory} \label{Linear_Theory}
\setcounter{equation}{0}

\bigskip
\bigskip
As it is mentioned in the introduction, weak solutions to the non-stationary Stokes system
\begin{equation}
\left\{ \ \gathered
\partial_t u-\Delta u+\nabla p= 0  \\
 \div u=0 \\
\endgathered \right. \quad \mbox{in} \quad Q
\label{Stokes_system}
\end{equation}
 locally are not necessary smooth. A simple example of a non-smooth
 solution to \eqref{Stokes_system} is as follows:
$$
u(x,t) = \ph(t)\nabla h(x), \qquad p(x,t)= -\ph'(t) h(x),
$$
where $h$ is a scalar harmonic function in spatial variables and $\ph$ is an arbitrary function of $t$ having limited smoothness. The same effect takes place in non-linear case and has been pointed out
  by J.~Serrin in  \cite{Serrin}.

Nevertheless,  Stokes system \eqref{Stokes_system} has the property of infinite smoothing of  weak solutions   with respect to spatial variables in internal points of a domain:

\begin{theorem} \label{T1}
Assume $s\in ]1, +\infty[$, $l\in ]1,2[$,  and  $u\in W^{1,0}_{s,l}(Q)$, $p\in L_{s,l}(Q)$ satisfy \eqref{Stokes_system} in $Q$ in the sense of distributions. Then  for any $k=0, 1, \ldots$, we have   $\nabla^k u \in C^{ \al, \frac \al 2}(\bar Q(\frac 12))$ with $\al= 2-\frac 2l$.
\end{theorem}

Surprisingly, the analog of Theorem \ref{T1} is not valid if we consider  weak solutions to the Stokes system near the boundary
\begin{equation}
\left\{
\gathered
 \ \gathered
\partial_t u-\Delta u+\nabla p= f  \\
 \div u=0 \\
 \endgathered  \quad \mbox{in} \quad Q^+
 \\
 u|_{x_3=0}=0. \qquad \qquad
\endgathered \right.
\label{Stokes_system_boundary}
\end{equation}
Actually, in contrast to the internal case, the first spatial gradient of a weak solution to system \eqref{Stokes_system_boundary} is not necessary bounded  up to  the boundary, i.e. there exist functions  $u\in W^{1,0}_2(Q^+)$, $p\in L_{\frac 32}(Q^+)$ which satisfy  system \eqref{Stokes_system_boundary} with $f\equiv 0$
in $Q^+$ in the sense of distributions, $u$ satisfy the boundary condition in the sense of traces but
\begin{equation}
\nabla u \not \in L_\infty(Q^+(r))
\label{Not_bounded}
\end{equation}
for any $0<r\leq 1/2$. The first counterexample  of this kind  has been constructed by Kang in  \cite{Kang 2005}. Later Seregin and Sverak in \cite{SS_ZNS385} simplified his construction significantly. Here,  we explain the counter-example, following to \cite{SS_ZNS385}:

\begin{counterexample} \label{T2}
Assume $\ph(t)$ is an arbitrary function of $t$ variable  and let $h: \Bbb R_+\times ]-~2, 0[\to \Bbb R$, $h=h(x,t)$ be a solution to the following initial boundary value problem for the  1D heat equation in a half-line:
$$
\left\{ \quad
\gathered
\frac{\cd h}{\cd t} - \frac{\cd^2 h}{\cd x^2} = \ph(t) \quad \mbox{in} \quad \Bbb R_+\times ]-2, 0[ \\
h|_{t=-2}=0, \qquad h|_{x=0}=0.
\endgathered
\right.
$$
Let $u: \Bbb R^3_+\times ]-2,0[\to \Bbb R^3$, $p: \Bbb R^3_+\times ]-2,0[\to \Bbb R$ be functions representing   the following  shear flow along $x_1$- axe:
$$
u(x,t): = (h(x_3, t), 0, 0), \qquad p(x,t):=- \ph(t)x_1
$$
Then functions $u$ and $p$ satisfy (formally) both the Stokes and the Navier-Stokes systems in $Q^+$. Moreover, if we assume  that $\al \in ]\frac 13, \frac 12[$ and take
$$
\begin{array}c
\ph(t)=\frac 1{|t|^{1-\al}}
\end{array}
$$
then  $u\in W^{1,0}_{2}(Q^+)$ and  $p\in L_{\frac 32}(Q^+)$. 
  Functions $u$ and $p$ satisfy equations \eqref{Stokes_system_boundary} with $f=0$ in the sense of distributions and the boundary condition  in the sense of traces as well. However, $\nabla u$ is unbounded in any neighborhood of the origin. In particular, \eqref{Not_bounded} holds.
\end{counterexample}

Example \ref{T2} shows that, in contrast to the internal case, the Stokes system does not possess the property of significant improvement of the regularity of weak solutions up to the boundary. This is a serious obstacle which makes  the  theory of boundary regularity for the Navier-Stokes equations  different from the analogues theory in the internal case. The natural question that arises is what is the optimal regularity of weak solutions to the Stokes system up to the boundary in the local set-up. A certain answer to this question has been given in \cite{Seregin_ZNS271}, where  H\" older continuity of the velocity field $u$ up to the flat part of the boundary has been established for {\it strong solutions} to \eqref{Stokes_system_boundary}. But before we switch to the discussion of this result let us introduce some terminology to explain the difference between weak and strong solutions to \eqref{Stokes_system_boundary}.

\begin{definition}
\label{Generalized_Solutions} Assume $1< s,l <+\infty$ and $f\in L_l(-1, 0; W^{-1}_s(B^+))$.
We say that functions $u$ and $p$ are a {\it weak solution} of system \eqref{Stokes_system_boundary}, if they belong to the spaces
$$
u\in W^{1,0}_{s,l}(Q^+), \qquad p\in L_{s,l}(Q^+),
$$
$u$ and $p$ satisfy \eqref{Stokes_system_boundary} in the sense of
distributions and $u$ satisfies the boundary condition
 in the sense of traces.
\end{definition}

Note that, for any weak solution $u$ and $p$ to system \eqref{Stokes_system_boundary},  $\cd_t u \in L_l(-1,0;W^{-1}_{s}(B^+))$  and the following estimate holds:
\begin{equation}
\gathered
\| \cd_t u\|_{L_l(-1, 0; W^{-1}_s(B^+))} \ \le \\ \le \  C~\Big( \| f\|_{L_{l}(-1, 0; W^{-1}_s(B^+))}+ \| u\|_{W^{1,0}_{s,l}(Q^+)}+ \| p\|_{L_{s,l}(Q^+)}\Big)
\endgathered
\label{Negative_time_derivative}
\end{equation}
Here $W^{-1}_s(B^+)$ is the space dual to  $\overset{\circ}{W}{^1_{s'}}(B^+)$
and
$$
\| u \|_{L_l(-1,0; W^{-1}_s(B^+))} = \Big(\int_{-1}^0 \| u(\cdot,t) \|^l_{W^{-1}_s(B^+)}~dt\Big)^{1/l}.
$$

\begin{definition}
\label{Strong_Solutions} Assume $1<s,l <+\infty$ and $f\in L_{s,l}(Q^+)$.
We say that  functions $u$ and $p$ are  a {\it strong solution} to \eqref{Stokes_system_boundary} if they are a weak solution to \eqref{Stokes_system_boundary} and
$$
u\in W^{2,1}_{s,l}(Q^+), \qquad p\in W^{1,0}_{s,l}(Q^+).
$$
\end{definition}


The idea of showing H\" older continuity of strong solutions proposed in \cite{Seregin_ZNS271} is as
follows. We first show that strong solutions satisfy the usual local version of
coercive estimate given by the following theorem:

\begin{theorem} \label{Theorem_Local_Estimate}
Suppose  $s$, $l \in ]1,\infty[$. For any $f\in
L_{s,l}(Q^+)$  and for any     strong  solution $u \in W^{2,1}_{s,l}(Q^+)$, $p \in W^{1,0}_{s,l}(Q^+)$
to  system \eqref{Stokes_system_boundary} in $Q^+$, the following local estimate holds:
\begin{equation}
\gathered
\| u \|_{W^{2,1}_{s,l} (Q^+(\frac 12))} + \| \nabla p \|_{L_{s,l}
(Q^+(\frac 12))}  \le \\ \le C \left(\| f \|_{L_{s,l}(Q^+)}+ \|
\nabla u \|_{L_{s,l} (Q^+)} + \| p\|_{L_{s,l} (Q^+)}
\right)
\endgathered
\label{Local_Estimate}
\end{equation}
with 
a positive constant $C$, depending only on  $s$, $l$.
\end{theorem}

In contrast to the interior case, this is not a trivial statement. The first proof
given in \cite{Seregin_ZNS271} has been based on duality arguments and inspired  by paper \cite{Solonnikov_1972}. In the present paper, we show
that it can be deduced from more general statement proved in \cite{FilShil}. But before explaining our approach let us briefly described the main idea of getting H\" older continuity. Unlike it has been done in the internal case, in the boundary case we can not get the result  by gaining more derivatives in space variables. In
fact, this is even impossible because of the counterexample in Example \ref{T2}. But what we can do is to gain more integrability in space and to apply certain bootstrap arguments. So, we show the following:

\begin{theorem} \label{Theorem3}
Suppose  $s$, $l$, $m \in ]1,\infty[$,  $m\ge s$. For any  $f\in
L_{m,l}(Q^+)$ and for any  strong solution  $u \in W^{2,1}_{s,l}(Q^+)$, $p \in W^{1,0}_{s,l}(Q^+)$ to \eqref{Stokes_system_boundary} in $Q^+$, we have 
 $$
 \begin{array}c
  u \in W^{2,1}_{m,l}(Q^+(\frac 12)), \qquad\nabla p \in
L_{m,l}(Q^+(\frac 12))
\end{array}
$$
and the following local estimate holds:
\begin{equation}
\gathered
\| u \|_{W^{2,1}_{m,l} (Q^+(\frac 12))} + \| \nabla p \|_{L_{m,l}
(Q^+(\frac 12))}  \le \\ \le C \left(\| f \|_{L_{m,l}(Q^+)}+ \|
\nabla u \|_{L_{s,l} (Q^+)} + \| p\|_{L_{s,l} (Q^+)}
\right)
\endgathered
\label{Local_Estimate2}
\end{equation}
with a positive constant $C$, depending only on  $s$, $l$, $m$.
\end{theorem}

If the exponent $m$ in Theorem \ref{Theorem3} is sufficiently large then the H\" older continuity of $u$ follows from the imbedding theorems for  the anisotropic Sobolev spaces. Namely, the following   result is true:

\begin{theorem} \label{Theorem4}
Suppose  $s$, $l$, $m \in ]1,\infty[$,  $m> \frac {3l}{2(l-1)} $. For any   $f\in
L_{m,l}(Q^+)$ and for any  strong solution  $u \in W^{2,1}_{s,l}(Q^+)$, $p \in W^{1,0}_{s,l}(Q^+)$ to  system \eqref{Stokes_system_boundary}, we have
 $ u \in C^{\beta, \frac \beta2}(\bar Q^+(\frac 12))$ with $\be = 2-\frac 3m -\frac 2l$ and the following local estimate holds:
\begin{equation}
\gathered
\| u \|_{ C^{\be, \frac \be 2}(\bar Q^+(\frac 12)) } \ \le \ C \left(\| f \|_{L_{m,l}(Q^+)}+ \|
\nabla u \|_{L_{s,l} (Q^+)} + \| p\|_{L_{s,l} (Q^+)}
\right),
\endgathered
\label{Local_Estimate3}
\end{equation}
where a positive constant $C$ depends only on  $s$, $l$, $m$.
\end{theorem}

The important part is to show that any weak solution to \eqref{Stokes_system_boundary} is  actually strong one at least locally (here we use the terminology introduced in our Definitions \ref{Generalized_Solutions} and \ref{Strong_Solutions}).
For the interior case,  the corresponding claim  is known (see details, for example, in \cite{Ladyzhenskaya_Seregin}). For the boundary case, the analogues statement has been proven in \cite{Seregin_ZNS370} relatively recently. Here, our
proof  follows to \cite{Shilkin_Vyalov}. It  is some modification of the approach in  \cite{Seregin_ZNS370}, which  can be
used in more general situations. So, the result is as follows:

\begin{theorem} \label{Theorem2}
Assume  $s$, $l \in ]1,\infty[$. Then, for any $f\in L_{s,l}(Q^+)$  and for any   weak solution $
u \in W^{1,0}_{s,l}(Q^+)$, $p \in L_{s,l}(Q^+)$  to \eqref{Stokes_system_boundary} in $Q^+$,  the following statemets
$$
\begin{array}c
u \in W^{2,1}_{s,l}(Q^+(\frac 12)),\qquad p \in W^{1,0}_{s,l}(Q^+(\frac 12))
\end{array}
$$
hold and  functions $u$ and $p$ are a strong solution to  system \eqref{Stokes_system_boundary} in the half-cylinder $Q^+(\frac 12)$.
\end{theorem}
As our example shows, regularity results up to the boundary described by the above statements are in a sense optimal.

Now we come to the detailed proofs of the results above.

\medskip
\noindent{\bf Proof of Theorem \ref{Theorem_Local_Estimate}:} The proof presented bellow is borrowed from  \cite{FilShil}. We reproduce it here for the sake of completeness.
Take arbitrary $\rho$, $r$ such that
$\frac12\le\rho<r\le\frac{9}{10}$.
Consider a cut-off function $\zeta\in C^\infty(\bar Q^+)$ such that
$$
\gathered
0\le \zeta\le 1 \quad\mbox{in} \quad Q^+,
\qquad\zeta\equiv 1 \quad\mbox{in} \quad Q^+(\rho),
\qquad \zeta\equiv 0 \quad\mbox{in} \quad Q^+\setminus Q^+(r), \\
 \| \nabla^k\zeta \|_{L_\infty(Q^+)}\le \frac {C}{(r-\rho)^{k}},\quad k=1,2,
\quad \| \cd_t\nabla^k \zeta \|_{L_\infty(Q^+)}\le \frac {C}{(r-\rho)^k}, \quad k=0,1.
\endgathered
$$
Let $u$ and $p$ be a strong solution to  system \eqref{Stokes_system_boundary}. Then functions
$v:=\zeta u$, $q:=\zeta p$  satisfy the following initial-boundary value problem
\begin{equation}
\left\{ \quad
\gathered
\gathered
\cd_t  v - \Delta  v +  \nabla q = \tilde f \\
\div v = g\\
\endgathered
 \qquad \text{in} \quad \Om \times ]-1,0[, \\
v|_{\cd \Om \times ]-1,0[} = 0, \quad v|_{t=-1} = 0.
\endgathered \right.
\label{Nonzero_div}
\end{equation}
where
\begin{equation}
\gathered
\tilde f = \zeta f + u (\cd_t \zeta- \Delta \zeta) - 2(\nabla u)\nabla \zeta+ p\nabla \zeta , \qquad
g = u\cdot \nabla \zeta
\endgathered
\label{New_RHS}
\end{equation}
and  $\Om $ is some smooth canonical domain which is diffeomorphic to a ball and satisfies the inclusions \
$B^+_{{9}/{10}}\subset \Om\subset B^+$.
Applying  Theorem 1.1 of \cite{FilShil}, we obtain the estimate
\begin{equation}
\gathered
\| v \|_{W^{2,1}_{s,l} (\tilde Q)} + \| \nabla
q\|_{L_{s,l}(\tilde Q)} \le \\ \le C \left(\|\tilde f\|_{L_{s,l}(\tilde Q)}+ \|  g  \|_{W^{1,0}_{s,l}(\tilde Q)} +
\| \cd_t g\|_{L_{s,l}(\tilde Q)}^{1/s }  \|
\cd_t g \|_{L_{l}(-1, 0; W^{-1}_s(\Om))}^{1/s'}\right)
\endgathered
\label{Stokes_Nonzero}
\end{equation}
where we denote $\tilde Q:=\Om \times ]-1,0[$.
Taking into account that
$\frac 1{r-\rho}\ge 1$ after routine computations   we obtain the estimate
\begin{equation}
\gathered
\| u \|_{W_{s,l}^{2,1} (Q^+(\rho))}^s
\le C\| f\|_{L_{s,l}(Q^+)}^s  \\
+ \frac C{(r-\rho)^{2s}} \Big(\| u \|_{W^{1,0}_{s,l}(Q^+)}^s
+ \|p\|_{L_{s,l}(Q^+)}^s+ \| \cd_t u \|_{L_l(-1,0; W^{-1}_s(B^+))}^{s}\Big) \\
+ \frac C{(r-\rho)^{2s}} \| \cd_t u \|_{L_{s,l}(Q^+(r))}\Big(
\| \cd_t u \|_{L_l(-1,0; W^{-1}_s(B^+))}^{s-1}+ \|  u \|_{L_{s,l}(Q^+)}^{s-1}\Big) .
\endgathered
\label{To_chto_nado}
\end{equation}
Estimating the last term in the right-hand side of \eqref{To_chto_nado}
via the Young inequality $ab\le \ep a^s+C_\ep b^{s'}$ we obtain
\begin{gather*}
\| u \|_{W_{s,l}^{2,1} (Q^+(\rho))}^s
\le  C\|  f\|^s_{L_{s,l}(Q^+)}+ \ep \| \cd_t u \|_{L_{s,l}(Q^+(r))}^s + \\
 +
\frac {C_\ep}{(r-\rho)^{2ss'}} \Big(\| u \|_{W^{1,0}_{s,l}(Q^+)}^s  + \|p\|_{L_{s,l}(Q^+)}^s+ \| \cd_t u \|_{L_l(-1,0; W^{-1}_s(B^+))}^{s}\Big) ,
\end{gather*}
where a constant $\ep>0$ can be chosen arbitrary small.
By virtue of \eqref{Negative_time_derivative}, we obtain
\begin{equation}
\gathered
\| u \|_{W_{s,l}^{2,1} (Q^+(\rho))}^s
\le  \ep \| \cd_t u \|_{L_{s,l}(Q^+(r))}^s
 +
\frac {C_\ep}{(r-\rho)^{2ss'}}
\Big(\| f\|_{L_{s,l}(Q^+)}^s + \| u \|_{W^{1,0}_{s,l}(Q^+)}^s + \|p\|_{L_{s,l}(Q^+)}^s \Big) .
\endgathered
\label{inequality this implies}
\end{equation}
Now, let us introduce the monotone function
$\Psi(\rho) := \| u \|_{W_{s,l}^{2,1} (Q^+(\rho))}^s$
and the constant
$$
A:=C_\ep\left(\| f\|_{L_{s,l}(Q^+)}^s + \| u \|_{W^{1,0}_{s,l}(Q^+)}^s  + \|p \|_{L_{s,l}(Q^+)}^s\right).
$$
The inequality \eqref{inequality this implies} implies that
\begin{equation}
\begin{array}c
\Psi (\rho) \le \ep\Psi(r)+\frac {A} {(r-\rho)^\al}, \qquad \forall ~\rho, \ r:
\quad R_1 \le \rho<r\le R_0,
\end{array}
\label{Giaquinta's lemma}
\end{equation}
for some  $\al>0 $ depending only on $s$, and for $R_1=\frac 12$,  $R_0=\frac 9{10}$.
Now we shall take an advantage of  the following lemma (which can be easily  proved by iterations if one take $r_k:=R_0-2^{-k}(R_0-R_1)$):

\begin{lemma}
Assume $\Psi$ is a nondecreasing bounded function which satisfies inequality \eqref{Giaquinta's lemma} for some $\al>0$, $A>0$, and $\ep\in ]0,2^{-\al}[$.
Then there exists a constant $B$ depending only on  $\ep$ and $\al$ such that
$$
\Psi(R_1)\le \frac {B\, A}{(R_0-R_1)^\al} .
$$
\label{Giaquinta_lemma}
\end{lemma}
Fixing $\ep = 2^{-3 ss'}$ in \eqref{inequality this implies}, applying
Lemma \ref{Giaquinta_lemma} to our function $\Psi$, and  evaluating  $\nabla p $ from  equations \eqref{Stokes_system_boundary} held a.e. in $Q^+$,
we derive the estimate
\eqref{Local_Estimate}. Theorem \ref{Theorem_Local_Estimate} is proved. $\square$



Theorems \ref{Theorem_Local_Estimate} together with results of \cite{FilShil} provides us the following proof of Theorem \ref{Theorem2}. 

\medskip
\noindent
{\bf Proof of Theorem \ref{Theorem2}.}  Let $\rho_m \to +0$ be an arbitrary sequence. Extend all functions $u$, $p$, $f$ from $Q^+$ to the set $B^+\times \Bbb R$ by zero. For any extended function $u$ denote by $u^m$ the mollification of the function $u$ with respect to $t$ variable:
$$
u^m(x,t) := (\om_{\rho_m}* u)(x,t) \equiv \int_{\Bbb R} \omega_{\rho_m} ( t - \tau ) u(x,\tau)\,d\tau,
$$
where $ \omega_{\rho}(t) = \frac1{\rho}\omega(t/\rho )$,
and $\omega \in C^\infty_0(-1, 1)$ is a smooth kernel normalized by the identity $\int_0^1 \om(t)dt = 1$.

As $u \in W^{1,0}_{s,l}(Q^+)$, $p \in L_{s,l}(Q^+)$, $f \in L_{s,l}(Q^+)$ we have
\begin{equation}
\gathered
u^m \to u \quad \text{ in } W^{1,0}_{s,l}(Q^+), \quad p^m \to p \quad \text{ in } L_{s,l}(Q^+),\\
f^m \to f \quad \text{ in } L_{s,l}(Q^+).
\endgathered
\label{wsl2}
\end{equation}
Let us fix arbitrary $\dl\in ]0,\frac 1{12}[$. Then for any $\rho_m< \dl$ and for any $\eta\in C^\infty(\bar Q^+)$
$$
\cd_t(\om_{\rho_m}* \eta) (x,t) = (\om_{\rho_m}* \cd_t \eta)(x,t), \quad \forall~x\in B^+, \ t\in ]-1+\dl, -\dl[.
$$
A weak solution $u$ and $p$ to  system \eqref{Stokes_system_boundary} obeys the integral identity
$$
\gathered
 -  \ \int\limits_{Q^+} u\cdot (\cd_t \eta +  \Delta \eta)~dxdt \ = \  \int\limits_{Q^+} (f\cdot \eta + p \div \eta)~dxdt
\endgathered
$$
which holds for all $\eta\in C^\infty(\bar Q^+)$ satisfying conditions
\begin{equation}
\gathered
\eta|_{\cd B^+\times ]-1,0[}=0, \qquad \nabla \eta |_{\cd' B^+\times ]-1,0[}=0, \qquad
\eta|_{B^+\times (]-1, -1+\dl[\cup ]-\dl, 0[)} =0,
\endgathered
\label{eta}
\end{equation}
where $\cd'B^+:=\{ x\in \Bbb R^n: |x|=1, x_n>0\}$.
Take 
the test function
$\eta= \om_{\rho_m}* \tilde \eta$, where $\tilde \eta\in C^\infty(\bar Q^+)$ is an arbitrary function  satisfying \eqref{eta}. Using properties of convolution,
we find the identity
\begin{equation}
\gathered
-~ \int\limits_{Q^+} u^m\cdot (\cd_t \tilde \eta +   \Delta \tilde \eta)~dxdt \ = \  \int\limits_{Q^+} (f^m\cdot  \tilde \eta + p^m \div  \tilde \eta)~dxdt
\endgathered
\label{Integral_Identity-1}
\end{equation}
which holds for all $\tilde \eta\in C^\infty(\bar Q^+)$  satisfying  \eqref{eta}.

Let $\zeta \in C^\infty(\bar{Q}^+)$ be a cut-of function vanishing in $ Q^+ \setminus Q^+(\frac 56)$ and such that $\zeta \equiv 1$ in $Q^+(\frac 23)$. Denote
$v^m := \zeta u^m$, $q^m := \zeta p^m$. Then from \eqref{Integral_Identity-1} we deduce that $v^m$ and $q^m$  obey the integral identity
$$
\gathered
-~ \int\limits_{B^+\times ]-1,-\dl[} v^m\cdot (\cd_t \eta +  \Delta \eta)~dxdt \ = \  \int\limits_{B^+\times ]-1,-\dl[} (\tilde f^m\cdot \eta + q^m \div  \eta)~dxdt
\endgathered
$$
for any $\eta\in C^\infty(\bar B^+\times [-1,-\dl])$ such that   $\eta|_{\cd B^+\times ]-1, -\dl [}=0$ and $\eta|_{B^+\times \{t=-\dl \}}=0$. Here $\tilde f^m$ and $g^m$
are determined by formulas \eqref{New_RHS} with $u$, $p$ and $f$ replaced by $u^m$, $p^m$ and $f^m$ respectively.

Assume $\Om\subset \Bbb R^3$ is a smooth domain such that $B^+(\frac 56)\subset \Om \subset B^+$ and denote $\tilde Q:=\Om \times ]-1,0[$.
As functions $g^m$ are  smooth with respect to $t$, we obtain   from Theorem 1.1 of \cite{FilShil} that, for any $m\in \Bbb N$, there exists a strong solution $\tilde v^m\in W^{2,1}_{s,l}(\tilde Q)$, $\tilde q^m\in W^{1,0}_{s,l}(\tilde Q)$ to the problem
\begin{equation}
\left\{ \quad
\gathered
\gathered
\cd_t \tilde v^m - \Delta \tilde v^m + \nabla \tilde q^m = \tilde f^m \\
\div  \tilde v^m = g^m\\
\endgathered
 \quad \text{in} \quad \tilde Q, \\
\tilde v^m|_{\cd\Om \times ]-1, 0[} = 0, \qquad \tilde v^m|_{t=-1}=0.
\endgathered \right.
\label{wsl3}
\end{equation}
Note that  $\zeta\equiv 1$ in $Q^+(\frac 23)$ and hence $g^m\equiv 0$ in $Q^+(\frac 23)$. So, functions $\tilde v^m$ and $\tilde q^m$ satisfy all assumptions of Theorem \ref{Theorem_Local_Estimate} in $Q^+(\frac 23)$ and, hence, by its obvious modification, 
we have the estimate
\begin{equation}
\gathered
\| \tilde{v}^m \|_{W^{2,1}_{s,l}(Q^+(\frac 12))} + \| \nabla \tilde{q}^m \|_{L_{s,l}(Q^+(\frac 12))} \ \le \\
\le \ C~ \left( \| \tilde f^m \|_{L_{s,l}(Q^+(\frac 23))} + \| \tilde v^m \|_{W^{1,0}_{s,l}(Q^+(\frac 23))} + \| \tilde q^m\|_{L_{s,l}(Q^+(\frac 23))} \right)
\endgathered
\label{wsl4}
\end{equation}
where a constant $C$  depends neither on $m$ nor on $\dl$.

 Since every strong solution of the  Stokes system is a weak one,  $\tilde v^m$ and $\tilde q^m$ satisfy the integral identity
$$
\gathered
-~ \int\limits_{\tilde Q}\tilde v^m\cdot (\cd_t \eta + \Delta \eta)~dxdt \ = \  \int\limits_{\tilde Q} (\tilde f^m\cdot \eta + \tilde q^m \div  \eta)~dxdt
\endgathered
$$
for all $\eta\in C^\infty(\overline{ \tilde Q})$ such that   $\eta|_{\cd \Om\times ]-1, 0 [}=0$ and $\eta|_{\Om \times \{t=0 \}}=0$. Hence the differences $w^m:=v^m-\tilde v^m$, $\pi^m:=q^m-\tilde q^m$  are a weak solution to the Stokes system   in $\Om\times ]-1, -\dl[$, satisfying the  identities
\begin{equation}
\gathered
\div  w^m =0 \quad \mbox{a.e. in}\quad \Om\times ]-1,-\dl[,
\\
-~ \int\limits_{\Om\times ]-1,-\dl[} w^m\cdot (\cd_t \eta +  \Delta \eta)~dxdt \ = \  \int\limits_{\Om\times ]-1,-\dl[} \pi^m \div \eta~dxdt,
\endgathered
\label{Integral_Identity-2}
\end{equation}
for any $\eta\in W^{2,1}_{s', l'}(\Om\times ]-1,-\dl[)$ such that   $\eta|_{\cd \Om\times ]-1, -\dl [}=0$ and $\eta|_{\Om\times \{t=-\dl \}}=0$.
Denote $\ka=\min\{s,l\}>1$. As $u^m$, $\tilde u^m\in L_{s,l}(\tilde Q)$ and $q^m$, $\tilde q^m\in L_{s,l}(\tilde Q)$ we have $w^m=v^m-\tilde v^m\in L_\ka(\tilde Q)$ and $\pi^m =q^m-\tilde q^m \in L_\ka(\tilde Q^+)$. Hence $|w^m|^{\ka-2}w^m\in L_{\ka'}(\tilde Q)$, and using results of \cite{Solonnikov_1972} we can find  functions $\eta\in  W^{2,1}_{\ka'}(\Om \times ]-1,-\dl[)$ and $\kappa \in W^{1,0}_{\ka'}(\Om \times ]-1,-\dl[)$ such that
$$
\left\{\quad
\gathered
\gathered
\cd_t \eta + \Delta \eta +\nabla \kappa = |w^m|^{\ka-2} w^m,\\
\div \eta = 0,
\endgathered  \quad \text{in} \quad \Om\times ]-1, -\dl[, \\
\eta|_{\cd \Om\times ]-1,\-\dl[}=0, \qquad \eta|_{t=-\dl}= 0.
\endgathered
\right.
$$
Substituting this $\eta$ as a test function into identity \eqref{Integral_Identity-2} we obtain $w^m =0$ in $\Om\times ]-1, -\dl[$. Hence
$v^m =\tilde v^m\in W^{2,1}_{s,l}(\Om\times ]-1, -\dl[)$. From \eqref{Integral_Identity-2} we obtain
\begin{equation}
\int\limits_{\Om\times ]-1,-\dl[} \pi^m \div \eta~dxdt \ = \ 0, \quad \forall ~\eta\in L_{l'}(-1,-\dl; \overset{\circ}{W}{^1_{s'}}(\Om)).
\label{Pressure}
\end{equation}
Correcting, if necessary, the function $\tilde q^m$ by a constant, we can assume that $\int\limits_\Om \pi^m~dx =0$ for a.e. $t\in ]-1,-\dl[$.
As $\pi^m \in L_\ka(\Om)$ for a.e. $t\in ]-1,-\dl[$, we have
$|\pi^m|^{\ka-2}\pi^m \in L_{\ka'}(\Om)$ for a.e. $t\in ]-1,-\dl[$. Using results of \cite{Bogovskii}  for a.e. $t$ we can find  $\eta(\cdot,t) \in  \overset{\circ}{W}{^1_{\ka'}}(\Om)$
such that
$$
\left\{ \quad
\gathered
\div \eta = |\pi^m|^{\ka-2}\pi^m -(|\pi^m|^{\ka-2}\pi^m)_\Om , \quad \mbox{a.e. } t\in ]-1, -\dl[, \\
\| \eta\|_{W^1_{\ka'}(\Om)} \le C\|\pi^m\|_{L_{\ka}(\Om)}^{\ka-1}.
\endgathered \right.
$$
From the latter estimate it follows that $\eta \in L_{\ka'}(-1,-\dl; \overset{\circ}{W}{^1_{\ka'}}(\Om))\subset L_{l'}(-1,-\dl; \overset{\circ}{W}{^1_{s'}}(\Om))$.
Substituting this $\eta$ into identity \eqref{Pressure},  we obtain
$\pi^m = 0$. This implies $q^m = \tilde q^m +const$ and we obtain the inclusion $q^m\in W^{1,0}_{s,l}(\Om\times ]-1, -\dl[)$. Moreover, from \eqref{wsl4},
we find
$$
\gathered
\| {v}^m \|_{W^{2,1}_{s,l}(B^+(\frac 12)\times ]-\frac {1}{4}, -\dl[)} + \| \nabla {q}^m \|_{L_{s,l}(B^+(\frac 12)\times ]-\frac {1}{4},-\dl[)} \ \le \\
\le \ C~ \left( \| \tilde f^m \|_{L_{s,l}(Q^+(\frac 23))} + \| v^m \|_{W^{1,0}_{s,l}(Q^+(\frac 23))} + \| q^m -b \|_{L_{s,l}(Q^+(\frac 23))} \right)
\endgathered
$$
where $C$ is independent on $m$ and $\dl$. Using identities $v^m=\zeta u^m$, $q^m=\zeta p^m$, $\zeta\equiv 1$ on $Q^+(\frac 23)$ and  expression \eqref{New_RHS} for $\tilde f^m$, we arrive at the estimate
$$
\gathered
\| u^m \|_{W^{2,1}_{s,l}(B^+(\frac 12)\times ]-\frac {1}{4}, -\dl[)} + \| \nabla p^m \|_{L_{s,l}(B^+(\frac 12)\times ]-\frac {1}{4},-\dl[)} \ \le \\
\le \ C~ \left( \| f^m \|_{L_{s,l}(Q^+(\frac 23))} + \| u^m \|_{W^{1,0}_{s,l}(Q^+(\frac 23))} + \| p^m  \|_{L_{s,l}(Q^+(\frac 23))} \right).
\endgathered
$$
Making use of \eqref{wsl2} we conclude that
$$
\begin{array}c
u\in W^{2,1}_{s,l}\left(B^+(\frac 12)\times ]-\frac {1}{4}, -\dl[\right), \quad p\in W^{1,0}_{s,l}\left(B^+(\frac 12)\times ]-\frac {1}{4},-\dl[\right),
\end{array}
$$
and the estimate
$$
\gathered
\| u \|_{W^{2,1}_{s,l}(B^+(\frac 12)\times ]-\frac {1}{4}, -\dl[)} + \| \nabla p \|_{L_{s,l}(B^+(\frac 12)\times ]-\frac {1}{4},-\dl[)} \ \le \\
\le \ C~ \left( \| f \|_{L_{s,l}(Q^+(\frac 23))} + \| u \|_{W^{1,0}_{s,l}(Q^+(\frac 23))} + \| p  \|_{L_{s,l}(Q^+(\frac 23))} \right)
\endgathered
$$
holds for any $\dl\in ]0, \frac 1{12}[$ with $C$ independent on $\dl$. The last inequality provides the required properties of $u$ and $p$. Theorem \ref{Theorem2} is proved.
$\square$



\medskip
Now, we are able to prove Theorem \ref{Theorem3}:

\medskip
\noindent
{\bf Proof of Theorem \ref{Theorem3}.}
For any $k=0,1,\ldots$ denote $s_k=\frac{ns}{n-ks}$ if $n>ks$ and $\frac{ns}{n-ks}<m$ and $s_k=m$ otherwise.  Denote also $N = \min\{ k\in \Bbb N: s_k = m\}$ and $\rho_k = \frac 12 + \frac 1{2^{k+1}}$.

Using obvious modification of Theorem 2.2 and Theorem 2.5, 
we observe that if $u \in W^{1,0}_{s_k, l}(Q^+(\rho_k))$ and  $p\in  L_{s_k,l}(Q^+(\rho_k))$ is a weak solution of  problem \eqref{Stokes_system_boundary} in $Q^+(\rho_k)$, then $u\in W^{2,1}_{s_k, l}(Q^+(\rho_{k+1}))$ and $p\in W^{1,0}_{s_{k},l}(Q^+(\rho_{k+1}))$ and the following estimate holds:
\begin{equation}
\gathered
\| u \|_{W^{2,1}_{s_k, l}(Q^+(\rho_{k+1}))} + \| \nabla p \|_{L_{s_k, l}(Q^+(\rho_{k+1}))} \ \le \\ \le \
C~\Big( \| f\|_{L_{m ,l}(Q^+)}  + \| u \|_{W^{1,0}_{s_k, l}(Q^+(\rho_k))} + \| p  \|_{L_{s_k, l}(Q^+(\rho_k))} \Big).
\endgathered
\label{1}
\end{equation}
Moreover, due to the imbedding $W^1_{s_k}(B^+(\rho_{k+1}))\hookrightarrow L_{s_{k+1}}(B^+(\rho_{k+1}))$, we find the estimate
\begin{equation}
\gathered
\| u \|_{W^{1,0}_{s_{k+1}, l}(Q^+(\rho_{k+1}))} + \| p \|_{L_{s_{k+1}, l}(Q^+(\rho_{k+1}))} \ \le  \\ \le
\ C~\Big(\| u \|_{W^{2,1}_{s_k, l}(Q^+(\rho_{k+1}))} + \| p \|_{W^{1,0}_{s_k, l}(Q^+(\rho_{k+1}))}\Big).
\endgathered
\label{2}
\end{equation}
Iterating \eqref{1} and \eqref{2} from $k=0$ to $k=N$ we finally obtain the bound
$$
\gathered
\| u \|_{W^{2,1}_{s_N, l}(Q^+(\frac 12))} + \|\nabla p \|_{L_{s_N, l}(Q^+(\frac 12))}  \ \le \\ \le \ C^N~
\Big( \| f\|_{L_{m ,l}(Q^+)}  + \| u \|_{W^{1,0}_{s_0, l}(Q^+)} + \| p  \|_{L_{s_0, l}(Q^+)} \Big).
\endgathered
$$
This estimate is equivalent to \eqref{Local_Estimate2}.
 Theorem \ref{Theorem3} is proved. \ $\square$



\medskip
\noindent
{\bf Proof of Theorem \ref{Theorem4}.} Theorem \ref{Theorem4} follows from Theorem \ref{Theorem3} and the following imbedding theorem for anisotropic Sobolev spaces (see \cite{Besov}):
$$
\gathered
 \begin{array}c
W^{2,1}_{m,l}(Q^+(\frac 12) )  \hookrightarrow C^{\be,\frac
\be 2}(\bar Q^+(\frac 12)), \qquad \mbox{if}\quad m> \frac {3l}{2(l-1)} \quad\mbox{and}\quad \be=2-\frac 3m-\frac 2l \end{array}
\\
\begin{array}c
\| u \|_{C^{\be,\frac \be 2}(\bar Q^+(\frac 12))}\le C \| u \|_{ W^{2,1}_{m,l}(Q^+(\frac 12 ))}, \qquad \forall~ u \in W^{2,1}_{m,l}(Q^+(\frac 12) ) .
\end{array}
\endgathered
$$
Theorem \ref{Theorem4} is proved. \ $\square$


\newpage

\section{Proof of the basic $\ep$--regularity condition} \label{Basic_Condition}
\setcounter{equation}{0}

\bigskip

\medskip
In this section we consider the Navier-Stokes system in a half-cylinder $Q^+$
\begin{equation}
\left\{ \quad
\gathered
\gathered
\cd_t u+  (u\cdot \nabla)u-\Delta u + \nabla p  = 0 \\
\div u = 0
\endgathered
\qquad \mbox{in} \quad Q^+
\\
u|_{x_3=0}=0 \quad\qquad \qquad
\endgathered
\right.
\label{NSE_boubdary}
\end{equation}
The aim of this section is to provide a proof  of the following theorem which is the boundary analogue of $\ep$-regularity condition \eqref{Basic_Condition_Internal}:

\begin{theorem}\label{Fixed_r} For any given $\al\in ]0, \frac 23[$ there exists a constant $\ep_*>0$ depending only on $\al$ such that for  any boundary suitable weak solution $u$ and $p$ in   $Q^+$  subject to the condition
\begin{equation}
\int\limits_{Q^+} \Big(~ |u|^3+|p|^{\frac 32}~\Big) dxdt \  < \ \ep_*,
\label{Fixed_r_condition}
\end{equation}
the velocity field $u\in C^{\al, \frac \al2}(\bar Q^+( \frac {1}2))$.
\end{theorem}

\noindent
We prove Theorem \ref{Fixed_r} following the method developed in \cite{Seregin_JMFM}. This method is based on the indirect approach in the regularity theory (see terminology, for example, in  \cite{Giaquinta}) and its crucial step   is  the decay estimate of Theorem \ref{Main_Lemma}. 
For a direct proof of partial regularity in the Navier-Stokes theory, we refer to \cite{Seregin_Aa}.

\begin{theorem}\label{Main_Lemma}
For any  $\theta\in ]0,\frac 12[$,  $\be\in ]0, \frac 23[$, there exists  a constant $\ep_0(\theta, \be)>0$  such that, for any boundary suitable weak solution  $u$ and $p$ to system \eqref{NSE_boubdary} in $Q^+$, the following implication holds:
$$
\mbox{if}\qquad Y_1(u,p)<\ep_0 \qquad \mbox{then} \qquad
Y_\theta(u,p) \le C_*~\theta^{\be}~Y_1(u,p).
$$
Here $C_*>0$ is some  absolute constant.
\end{theorem}

\noindent
Here we denote
$$
Y_\theta (u,p) \ := \  \ \Big(\pint\limits_{Q^+(\theta)} ~|u|^3~dxdt\ \Big)^{1/3} \ + \
\theta ~\Big( \pint\limits_{Q^+(\theta)} ~|p-[p]_{B^+(\theta)}|^{3/2}~dxdt\ \Big)^{2/3},
$$
where for any $f\in L_1(Q^+(\theta))$ we denote
$$
  \pint_{Q^+(\theta)} f(x, t)~dxdt = \frac 1{|Q^+(\theta)|}\int\limits_{Q^+(\theta)} f(x, t)~dxdt, \quad [f]_{B^+(\theta)} =  \frac 1{|B^+(\theta)|}\int\limits_{B^+(\theta)} f(x, t)~dx
$$

\medskip
\noindent
{\bf Proof of Theorem \ref{Main_Lemma}:}   Arguing by contradiction, we assume there
is a number $\theta\in ]0,\frac 12[$,  the sequence $\ep_h\to 0$
and functions $u^h$ and $p^h$ which are the  boundary suitable
weak solutions in the sense of Definition \ref{Definition_2} satisfying relations
$$
 Y_1(u^h,p^h) =
\ep_h \to 0, \qquad Y_\theta (u^h,p^h) \ge C_* \theta^{\be}
\ep_h.
$$
 We introduce  new functions
$$
v^h=\frac 1{\ep_h} u^h, \qquad q^h= \frac 1{\ep_h} (p^h-
[p^h]_{B^+}).
$$
They meet relations
\begin{equation}
Y_1(v^h,q^h) = 1, \qquad Y_\theta (v^h,q^h) \ge C_* \theta^{\be},
\label{Y_1=1}
\end{equation}
 as well as the system
\begin{equation}
\left\{ \quad \gathered
\gathered
\cd_t v^h + \ep_h \div  (v^h \otimes
v^h) - \Delta v^h + \nabla q^h = 0 \\
\div v^h = 0
\endgathered
\qquad\mbox{in } \  Q^+ ,
\\
v^h|_{x_3=0}=0, \qquad \qquad \qquad
\endgathered \right.
\label{NSE_u}
\end{equation}
which holds in the sense of distributions and the boundary condition is understood in the sense of traces.
Moreover, functions $v^h$ and $q^h$ satisfy the local energy inequality
\begin{equation}
\gathered
\int\limits_{B^+} \zeta(x,t)|v^h (x,t)|^2~dx + 2
\int\limits_{-1}^t \int\limits_{B^+} \zeta |\nabla v^h |^2~dxdt \ \le \\
\le    \ \int\limits_{-1}^t\int\limits_{B^+} \left\{~ |v^h
|^2\left(\cd_t \zeta + \Delta \zeta\right) +
v^h\cdot \nabla \zeta\left(\ep_h |v^h|^2+2q^h \right)
~ \right\}~dxdt
\endgathered
\label{LEI_u}
\end{equation}
 for a.e. $t\in ]-1,0[$ and all nonnegative $\zeta\in
C^\infty(\bar Q)$ vanishing near $\cd' Q$.

From (\ref{Y_1=1}), we derive the estimate
\begin{equation}
\| v^h \|_{L_3(Q^+)} + \|q^h\|_{L_{\frac 32} (Q^+) } \le C.
\label{weak_norms}
\end{equation}
 Picking up a cut-off function $\zeta$ so that $\zeta\equiv 1$ on $Q^+(\frac 34)$ and  taking into account (\ref{weak_norms}), we
find
\begin{equation}
\| v^h\|_{L_{2,\infty}(Q^+(\frac 34))} + \|
v^h\|_{W^{1,0}_2(Q^+(\frac 34))} \le C.
\label{energy_norms}
\end{equation}
The known multiplicative inequality allows us to conclude that
\begin{equation}
\| v^h\|_{L_{\frac {10}3}(Q^+(\frac 34))}\le C. \label{L_10/3}
\end{equation}
Another bound easily follows from (\ref{NSE_u}) and has the form 
\begin{equation}
\begin{array}c
\| \cd_t v^h \|_{L_{\frac 32}(-(\frac 34)^2,0;
W^{-1}_{\frac 32}(B^+(\frac 34)))}\le C. \label{time_derivative_u}
\end{array}
\end{equation}
Estimates \eqref{weak_norms}, \eqref{energy_norms} provide the
existence of
 subsequences $\{ v^h\}$ and $\{q^h\}$ with the following properties 
 \begin{equation}
 \begin{array}c
 v^h \rightharpoonup v^0 \quad \mbox{ in }\quad W^{1,0}_2(Q^+(\frac
 34))\cap L_3(Q^+),
\end{array}
 \label{weak_conv_u}
 \end{equation}
 \begin{equation}
 \begin{array}c
q^h \rightharpoonup q^0 \quad \mbox{ in }\quad L_{\frac 32}(Q^+).
 \end{array}
 \label{weak_conv_p}
\end{equation}
 Routine compactness arguments imply
 \begin{equation}
 \begin{array}c
   v^h \to v^0 \quad \mbox{ in }\quad L_3(Q^+(\frac
 34)),
 \end{array}
 \label{strong_conv_u}
\end{equation}
  Convergence (\ref{weak_conv_u}) --- (\ref{weak_conv_p}) allow us to
pass to the limit in  equations (\ref{NSE_u}) (if we take these
equations in the weak form). So,  $v^0$ and $q^0$
is a weak solution to the system
\begin{equation}
\left\{ \quad
\gathered
\gathered
\cd_t v^0 - \Delta v^0 + \nabla q^0 = 0
\\
\div  v^0 = 0
\endgathered
\quad\mbox{ in } \begin{array}c Q^+(\frac 34), \end{array} \\
v^0|_{x_3=0}=0.
\endgathered \right.
\label{NSE_u0}
\end{equation}
 Moreover, from the second
relation in (\ref{Y_1=1}),  we deduce the estimate
\begin{equation}
\liminf\limits_{h\to\infty} Y_\theta (v^h,q^h) \ge C_*
\theta^{\be}. \label{liminf}
\end{equation}
On the other hand, below we will show that
\begin{equation}
\limsup\limits_{h\to\infty} Y_\theta (v^h,q^h) \le C_{**}
\theta^{\be}, \label{limsup}
\end{equation}
with a  constant $C_{**}>0$. Taking  in
(\ref{liminf}) a constant $C_*>C_{**}$  we arrive at a contradiction
between (\ref{liminf}) and (\ref{limsup}). This  will complete our
proof Theorem \ref{Main_Lemma}.

\medskip
To prove
 (\ref{limsup}),  we split $Y_\theta(v^h, q^h)$ onto two parts:
$$
Y_\theta (v^h,q^h) = Y^1_\theta (v^h) + Y^2_\theta (q^h),
$$
where
$$
\gathered
Y^1_\theta (v^h)\equiv \Big( \pint\limits_{Q^+(\theta)} |v^h|^3
~dxdt \Big)^{\frac 13}, \qquad
 Y^2_\theta (q^h)\equiv \theta \Big( \pint\limits_{Q^+(\theta)}
|q^h |^{\frac 32}  ~dxdt \Big)^{\frac 23}
\endgathered
$$
As $\theta\in ]0,\frac 12[$,  strong convergence \eqref{strong_conv_u} gives us the following 
\begin{equation}
\lim\limits_{h\to \infty} Y^1_\theta(v^h) = Y^1_\theta(v^0).
\label{Comb1}
\end{equation}
 Since $v^0\in W^{1,0}_2(Q^+(\frac 34))$, $q^0\in L_{\frac 32}(Q^+(\frac 34))$ are a weak solution to the  Stokes system
(\ref{NSE_u0}) in $Q^+(\frac 34)$ (in the sense of Definition \ref{Generalized_Solutions}),  one can apply  Theorems \ref{Theorem2} and \ref{Theorem4}
find that $v^0\in C^{\be,\frac
\be2}(\bar Q^+(\frac 12))$ with the estimate
\begin{equation}
\| v^0 \|_{C^{\be,\frac \be2}(\bar
Q^+(\frac 12))}  \le C \left(\|
\nabla v^0 \|_{L_{\frac 32} (Q^+(\frac 34))} + \| q^0 \|_{L_{\frac
32}(Q^+)} \right). \label{est_u0}
\end{equation}
Thanks to estimates (\ref{weak_norms}), (\ref{energy_norms})
and the lower semicontinuity of the corresponding norms with
respect to  weak convergence (\ref{weak_conv_u}),
(\ref{weak_conv_p}), the right-hand side of (\ref{est_u0}) can be estimated by
some absolute constant $C$.
 As $v^0|_{x_3=0}=0$, relation
(\ref{est_u0}) yeilds
the estimate
\begin{equation}
\begin{array}c
Y^1_\theta(v^0)
 \le  C \theta^{\be} \|v^0\|_{C^{\be,\frac \be2}(\bar
Q^+(\frac 12))} \le C_0 \theta^{\be}.
\end{array}
\label{Comb2}
\end{equation}
 Combining (\ref{Comb1}) and (\ref{Comb2}),  we find the estimate
 \begin{equation}
\limsup\limits_{h\to\infty} Y_\theta (v^h,q^h) \le C_0
\theta^{\be} + \limsup\limits_{h\to\infty} Y^2_\theta (q^h)
\label{Comb3}
 \end{equation}
So, to show (\ref{limsup}), we need to estimate the second
term in the right-hand side of (\ref{Comb3}). For this purpose, let us define
$$
f^h=
-\ep_h \div (v^h \otimes  v^h) \quad \mbox{in}\quad  \begin{array}c Q^+(\frac 34). \end{array}
$$
 From  H\" older inequality, it follows that  $f^h \in L_{\frac 98, \frac 32}(Q^+(\frac 34))$ and using \eqref{energy_norms} we obtain
\begin{equation}
\| f^h \|_{L_{\frac 98, \frac 32}(Q^+(\frac 34))} \ \le  \ C~\ep_h~ \| v^h \|_{L_{2,\infty}(Q^+(\frac 34))}^{\frac 23}\| \nabla v^h \|_{L_{2}(Q^+(\frac 34))}^{\frac 43}
\to 0 \quad \mbox{as}\quad h\to \infty.
\label{f_tens_to_zero}
\end{equation}
Assume $\Om\subset \Bbb R^3$ is a smooth canonical domain such that $B^+(\frac 58)\subset \Om \subset B^+(\frac 34)$ and denote $\tilde Q:=\Om \times ]-\frac {9}{16},0[$. It is known, see, for example, \cite{Solonnikov_UMN}, that there exist 
a unique pair of functions $\hat v^h\in W^{2,1}_{\frac 98, \frac 32}(\tilde Q)$,  $q^h_1\in W^{1,0}_{\frac 98, \frac 32}(\tilde Q)$, $[q^h_1]_{\Om}=0$ a.e. $t\in ]-\frac 9{16},0[$,  which obey the following initial-boundary value problem
$$
\left\{ \quad
\gathered
\gathered
\cd_t \hat v^h   - \Delta \hat v^h +
 \nabla q^h_1 =f^h \\ \div   \hat v^h=0
  \endgathered
\quad\mbox{ in }\quad \tilde Q, \\
\hat  v^h|_{t=-\frac 9{16}}=0, \qquad \hat
v^h|_{\cd \Om \times ]-\frac 9{16},0[}=0,
\endgathered \right.
$$
and is subject to the estimate
\begin{equation}
\begin{array}c
\| \hat  v^h \|_{W^{2,1}_{\frac 98, \frac 32}(\tilde Q)}
 + \|  q^h_1 \|_{W^{1,0}_{\frac 98,
\frac 32}(\tilde Q)} \ \le \   C \| f^h \|_{L_{\frac 98, \frac
32}(\tilde Q)}.
\end{array}
\label{14}
\end{equation}
From \eqref{f_tens_to_zero}, \eqref{14}, and from the imbedding $W^1_{\frac 98}(\Om)\hookrightarrow L_{\frac 32}(\Om)$, 
we can conclude
\begin{equation}
Y^2_\theta(q_1^h ) \le C\theta^{-2} \| \nabla q_1^h\|_{L_{\frac 98, \frac 32}(\tilde Q)}  \ \le \ C~\theta^{-2} \| f^h\|_{L_{\frac 98, \frac 32}(Q^+(\frac 34))}   \to 0 \quad \mbox{as}\quad h\to \infty.
\label{q_1}
\end{equation}

Now, consider the functions $ \tilde v^h \equiv  v^h- \hat v^h $, $
q^h_2 = q^h - q^h_1$. Note that $\tilde v^h \in W^{1,0}_{\frac 98, \frac 32}(\tilde Q)$, $q^h_2 \in L_{\frac 98, \frac 32}(\tilde Q)$ and, hence, $\tilde v^h$ and $q^h_2$ is a weak solution (in the sense of Definition \ref{Generalized_Solutions}) of the homogeneous Stokes system in $\tilde Q$
$$
\left\{ \quad \gathered
\gathered
\cd_t \tilde v^h  - \Delta \tilde v^h +
\nabla q^h_2 =0 \\  \div   \tilde v^h=0
\endgathered
\quad\mbox{
in }\quad \tilde Q, \\
 \tilde
v^h|_{x_3=0}=0.
\endgathered
\right.
$$
Let us take $m:=\frac {9}{2-3\be}$, $m\in ]\frac 92, +\infty[$.
By an obvious modification of Theorem   \ref{Theorem3} with $f\equiv 0$ we obtain inclusions $\tilde v^h \in W^{2,1}_{m, \frac 32}(Q^+(\frac 12))$, $q_2^h \in W^{2,1}_{m, \frac 32}(Q^+(\frac 12))$ and  the  estimate
$$
\| \nabla q^h_2 \|_{L_{m,\frac 32}(Q^+(\frac 12))} \le C
\left(\|  \tilde v^h \|_{W^{1,0}_{\frac 98, \frac 32}(\tilde Q)} +   \|  q^h_2
\|_{L_{\frac 98, \frac 32}(\tilde Q)}\right).
$$
The right-hand side  of the last inequality can be controlled by
majorants which do not depend on $h$. Indeed, as $\tilde v^h = v^h -\hat v^h$, from \eqref{energy_norms}, \eqref{14}, we find
$$
\| \nabla \tilde v^h \|_{L_{\frac 98, \frac 32}(Q^+(\frac 34))}
\le \| \nabla  v^h \|_{L_{2}(Q^+(\frac 34))} + \| \nabla \hat v^h
\|_{L_{\frac 98, \frac 32}(Q^+(\frac 34))} \le C.
$$
As $q^h_2=q^h - q^h_1$, using \eqref{weak_norms}, \eqref{f_tens_to_zero},  \eqref{14} we get
$$
\gathered
\|  q^h_2  \|_{L_{\frac 98,
\frac 32}(\tilde Q)}  \le C \left(\| q^h \|_{L_{\frac 32}(Q^+) } + \| q^h_1  \|_{L_{\frac 98, \frac 32}(\tilde Q)}\right) \le
 C\left(1+ \| f^h \|_{L_{\frac 98, \frac 32}(Q^+(\frac 34))}\right)\le C .
\endgathered
$$
 Hence, it has been shown that
 \begin{equation}
\| \nabla q^h_2 \|_{L_{m,\frac 32}(Q^+(\frac 12))} \le C .
\label{High_norm_of_pressure}
\end{equation}
Using  \eqref{High_norm_of_pressure} and the Poincar\'{e} and H\" older inequalities, we find
\begin{equation}
\begin{array}c
Y^2_\theta(q^h_2)  \le  C  \theta^{\be } \| \nabla
q^h_2 \|_{L_{m,\frac 32}(Q^+(\frac 12))} \le C_1 \theta^{\be}.
\end{array}
\label{q_2}
\end{equation}
Combining \eqref{q_2} with  \eqref{q_1}, we show
$$
\limsup\limits_{h\to \infty}Y^{2}_\theta(q^h) \le \limsup\limits_{h\to \infty}Y^{2}_\theta(q^h_1) +\limsup\limits_{h\to \infty}Y^{2}_\theta(q^h_2)  \le 0 + C_1\theta^\be= C_1\theta^\be.
$$
Hence from \eqref{Comb3} we obtain
(\ref{limsup})  with $C_{**}=C_0+C_1$
 which contradicts to (\ref{liminf}) if we take $C_*>C_{**}$.
Theorem \ref{Main_Lemma} is proved. \ $\square$

\begin{theorem} \label{Scaling_Theorem}
Let $C_*>1$ be the absolute constant defined by Theorem \ref{Main_Lemma}. Assume $\be\in ]0, \frac 23[$, $\theta\in ]0, \frac 12[$ are arbitrary. Denote by $\ep_0>0$, $\ep_0=\ep_0(\theta, \be )$ a constant determined by Theorem \ref{Main_Lemma}. Then, for any boundary suitable weak solution $u$ and $p$ to the Navier-Stokes system in $Q^+$ and for any $k=0,1,2,\ldots$,  the following is true:
\begin{equation}
\mbox{if}\quad Y_{\theta^k}(u,p)<\ep_0 \qquad \mbox{then} \qquad    Y_{\theta^{k+1}}(u,p) \le  C_* \theta^\be ~Y_{\theta^k}(u,p)
\label{Conclusion_Scaling_Theorem}
\end{equation}
\end{theorem}

\medskip
\noindent
{\bf Proof:} Define functions $u^{\theta^k}$ and $p^{\theta^k}$ by formulas
$$
\gathered
u^{\theta^k}(x,t) \ := \ \theta^k u(\theta^k x, \theta^{2k}t) \\
p^{\theta^k}(x,t) \ := \ \theta^{2k} p(\theta^k x, \theta^{2k}t)
\endgathered
\qquad (x,t)\in Q^+.
$$
As $u$ and $p$ is a boundary suitable weak solution of the Navier-Stokes equations in  $Q^+(\theta^k)$, the functions $u^{\theta^k}$ and $p^{\theta^k}$ are a boundary suitable weak solution of the Navier-Stokes system in  $Q^+$. Moreover,
$$
Y_1(u^{\theta^k}, p^{\theta^k}) \ = \ Y_{\theta^k}(u,p) \ < \ \ep_0,
$$
and, hence, by Theorem \ref{Main_Lemma}
$$
Y_\theta (u^{\theta^k}, p^{\theta^k})  \ \le \ C_* \theta^\al ~Y_1(u^{\theta^k}, p^{\theta^k}).
$$
From this inequality, the conclusion of the implication  \eqref{Conclusion_Scaling_Theorem} follows by change of variables. Theorem \ref{Scaling_Theorem} is proved. \ $\square$

\medskip

\begin{theorem} \label{Iterations}
Let $C_*>1$ be the absolute constant defined by Theorem \ref{Main_Lemma} and let $\al\in ]0,\frac 23[$, $\be\in ]\al, \frac 23[$ be arbitrary. Assume a number $\theta\in ]0, \frac 12[$ is fixed in such a way that
\begin{equation}
C_* \theta^{\be-\al} <1,
\label{Choice of theta}
\end{equation}
and let $\ep>0$, $\ep_0=\ep_0(\theta, \be)$ be  a constant determined by Theorem \ref{Main_Lemma}.
Then, for any boundary suitable weak solution $u$ and $p$ of the Navier-Stokes system in $Q^+$ and for any $k=0,1,2,\ldots$, the following is valid:
$$
\mbox{if}\quad Y_1(u,p)<\ep_0 \quad \mbox{then for any} \quad k=0,1,2,\ldots \quad
\left\{ \quad
\gathered
Y_{\theta^k}(u,p) <\ep_0 \\
Y_{\theta^{k+1}}(u,p) < \theta^{\al (k+1)} Y_1(u,p)
\endgathered
\right.
$$
\end{theorem}

\medskip
\noindent
{\bf Proof:} The proof follows easily from Theorem \ref{Scaling_Theorem} by induction in $k$.  \ $\square$

\medskip

\begin{theorem} \label{Shifts}
 Assume $\al\in ]0, \frac 23[$ is arbitrary and take $\be = \frac 12(\al+\frac 23)$. Let us fix  $\theta\in ]0, \frac 12[$, $\theta=\theta(\al)$ so that \eqref{Choice of theta} holds and let $\ep_0>0$, $\ep_0=\ep_0(\theta, \be)$ be  a constant determined by Theorem \ref{Main_Lemma}. Denote $\ep_*':= \ep_0(\theta, \be)$ and note that $\ep_*'>0$ actually depends only on $\al$.
Then for any $z_0=(x_0, t_0)$, $x_0\in \cd \Bbb R^3_+$, and for any  boundary suitable weak solution $u$ and $p$ of the Navier-Stokes equations in  $Q^+(z_0,R)$, the following is valid:
$$
\gathered
\mbox{if}\qquad R ~Y_{z_0, R}(u,p)<\ep_*' \\ \mbox{then for any} \quad 0<r<R \qquad
Y_{z_0, r}(u,p)  \ \le  \ C(\al)\left( \frac rR\right)^{\al} Y_{z_0, R}(u,p).
\endgathered
$$
Here,
$$
Y_{z_0,R} (u,p) \ := \  \ \Big(\pint\limits_{Q^+(z_0,R)} ~|u|^3~dxdt\ \Big)^{\frac 13} \ + \
R ~\Big( \pint\limits_{Q^+(z_0, R)} ~|p-[p]_{B^+(x_0,R)}|^{\frac 32}~dxdt\ \Big)^{\frac 23}
$$
\end{theorem}

\medskip
\noindent
{\bf Proof:} The proof follows easily from Theorem \ref{Iterations} by scaling transformation \eqref{Scale_Invariance}, if we fix $k\in \Bbb N\cup\{0\}$ in such a way that $\theta^{k+1}R \le r< \theta^k R$. $\square$

\medskip
\noindent
{\bf Proof of Theorem \ref{Fixed_r}:} Let
$$
\gathered
\bar Y_{z_0, R}(u,p) \ :=  \
 \Big(\pint\limits_{Q(z_0,R)\cap Q^+} ~|u-(u)_{Q(z_0,R)\cap Q^+}|^3~dxdt\ \Big)^{\frac 13} \ + \\ + \
R ~\Big( \pint\limits_{Q(z_0, R)\cap Q^+} ~|p-[p]_{B(x_0,R)\cap B^+}|^{\frac 32}~dxdt\ \Big)^{\frac 23},
\endgathered
$$
where $(u)_{Q(z_0,R)\cap Q^+}$ denotes the space-time average of $u$ over the set $Q(z_0,R)\cap Q^+$, and $[p]_{B(x_0,R)\cap B^+}$ denotes the spatial average of $p$ over the set $B(x_0,R)\cap B^+$.

For any  $x_0\in \Bbb R^3_+$, \ $x_0=(x_1^0, x_2^0, x_3^0)$,  \ $z_0=(x_0, t_0)$, \ we denote by  $x_0'$  the point with coordinates  $(x_1^0, x_2^0, 0)$, $z_0':=(x_0', t_0)$, and $d(x_0):=  \dist \{ x_0, \cd \Bbb R^3_+\} = x^0_3 $.

From the internal regularity theory (see, for example, \cite{Seregin_UMN}) we know there is a  constant $\ep_*''>0$ depending only on $\al$ such that
for any $z_0\in Q^+$, $R>0$ the following implication holds:
\begin{equation}
Q(z_0, R)\subset Q^+,  \  \   R~\bar Y_{z_0, R}(u,p)<\ep_*'' \quad  \Longrightarrow \quad \bar Y_{z_0,r}(u,p) \ \le \ c~\left(\frac rR\right)^\al ~\bar Y_{z_0, R}(u,p)
\label{Internal Decay}
\end{equation}
On the other hand, Theorem \ref{Shifts} reads that there is a constant $\ep_*'>0$ depending only on $\al$ such that
for any $z_0'\in \cd \Bbb R^3_+$, $R>0$, we have:
\begin{equation}
Q^+(z_0', R)\subset Q^+, \ \   R~ Y_{z_0', R}(u,p)<\ep_*' \ \ \Longrightarrow \ \ Y_{z_0',r}(u,p) \ \le \ c~\left(\frac rR\right)^\al ~ Y_{z_0', R}(u,p)
\label{Boundary Decay}
\end{equation}
Besides, it is easy to see that there is an absolute constant $c_1$ such that for any $z_0\in  B^+(\frac 12)$, \ $z_0' \in \bar B^+(\frac 12)\cap \cd \Bbb R^3_+$
\begin{equation}
\bar Y_{z_0, \frac 14}(u,p)\le c_1 Y_1(u,p), \qquad Y_{z_0', \frac 14}(u,p)\le c_1 Y_1(u,p)
\label{c1}
\end{equation}
Assume
 $z_0=(x_0, t_0)\in \bar Q^+(\frac 12)$  is arbitrary and   $0<r<\frac 18$. Denote $d:= d(x_0)$. There are  three possible cases:

\smallskip
\quad Case 1: \ $0 \le d < r < \frac 18$

\smallskip
\quad Case 2: \ $0<r \le d < \frac 18$

\smallskip
\quad Case 3: \  $\frac 18 \le d \le  \frac 12$

\noindent
In Case 1 we have $\bar Y_{z_0,r}(u,p) \ \le \  2~  Y_{z_0',2r}(u,p)$.
Let us fix $\ep_{1}:=\frac {\ep'_*}{c_1}$ where $c_1$ is fixed in \eqref{c1}. Then if we assume $Y_1(u,p)<\ep_{1}$ from \eqref{c1} we obtain
$Y_{z_0', \frac 14}(u,p) <\ep_*'$ and hence with the help of \eqref{Boundary Decay} we obtain
$$
\bar Y_{z_0,r}(u,p) \ \le \  2~  Y_{z_0',2r}(u,p) \  \le \ 2~c~\left(\frac {2r}{\frac 14}\right)^\al ~Y_{z_0',\frac 14}(u,p) \ \le \ c \ep_*' r^\al.
$$
In Case 2 we have $\bar Y_{z_0,d}(u,p) \ \le \  2~  Y_{z_0',2d}(u,p)$. If assume $Y_1(u,p)<\ep_{1}$, \eqref{c1} yields
$Y_{z_0', \frac 14}(u,p) <\ep_*'$ and hence with the help of \eqref{Boundary Decay}, \eqref{c1} we state that:
$$
\bar Y_{z_0,d}(u,p) \ \le \  2~  Y_{z_0',2d}(u,p) \  \le \ 2~c~\left(\frac {2d}{\frac 14}\right)^\al ~Y_{z_0',\frac 14}(u,p) \ \le \ c_2 Y_1(u,p) d^\al.
$$
Taking into account that $d<\frac 18$, we can conclude that $d~\bar Y_{z_0,d}(u,p) \le c_2 Y_1(u,p)$. Hence, if we fix $\ep_2:=\min \{ \frac {\ep'_*}{c_1}, \frac {\ep_*''}{c_2}\}$ and assume $Y_1(u, p)<\ep_2$, we can apply \eqref{Internal Decay} with $R=d$ and show
$$
\bar Y_{z_0,r}(u,p) \ \le \ c~\left(\frac rd\right)^\al ~\bar Y_{z_0, d}(u,p) \ \le \ c~\ep''_* ~r^\al
$$
In Case 3 we have $d~\bar Y_{z_0, d}(u,p)\le c_3Y_1(u,p)$ with an absolute constant $c_3>0$. Hence if we take $\ep_3:=\frac {\ep_*''}{c_3}$ and assume $Y_1(u,p)< \ep_3$ we observe that  $d~Y_{z_0, d}(u,p)< \ep_*''$ and hence we can apply \eqref{Internal Decay} with $R=d$. Taking into account $\frac 18\le d\le \frac 12 $, we find
$$
\bar Y_{z_0,r}(u,p) \ \le \ c~\left(\frac r{d}\right)^\al ~\bar Y_{z_0,d}(u,p) \le c ~\ep_*'' ~r^\al.
$$
Finally, if we fix $\ep_*:=\min\{ \ep_1, \ep_2, \ep_3\}$ and assume $Y_1(u,p)<\ep_*$ then in all cases we get the estimate
$$
\begin{array}c \forall~z_0\in \bar Q^+(\frac 12),  \quad \forall~r\in (0, \frac 18) \end{array} \qquad \bar Y_{z_0, r}(u,p)\le K r^\al,
$$
with some $K>0$ depending only on $\al$, $\ep_*'$ and $\ep_*''$. From this estimate, we deduce H\" older continuity of $u$ in the set $\bar Q^+(\frac 12)$ via Campanato criterion, see, for example, \cite{Campanato}. Theorem \ref{Fixed_r} is proved. \ $\square$

\newpage

\section{Further results} \label{Curved_Boundary}
\setcounter{equation}{0}

\medskip
\bigskip
In this section we discuss further results of $\ep$-regularity theory for the Navier-Stokes system.
We start with the Navier-Stokes equations in the neighborhood of a  point $x_0$ belonging to  the smooth curvilinear part of the  boundary $\cd\Om$ of a domain $\Om\subset \Bbb R^3$.
Namely, assume $x_0\in \cd \Om$, denote by $\Om(x_0, R)$ the intersection of some neighborhood of $x_0$ with $\Om$ and consider the system
\begin{equation}
\left\{ \quad \gathered \gathered \cd_t u - \Delta u + (u\cdot \nabla )u +\nabla p  \ = \  f \\ \div u = 0 \endgathered \qquad \mbox{in}\quad \Om(x_0, R)\times ]-R^2,0[.  \\
u|_{\cd\Om\times ]-R^2, 0[}=0 \qquad\qquad  \qquad\qquad  \qquad\qquad \quad
\endgathered\right.
\label{Stokes_system-1}
\end{equation}
Without loss of generality we can assume that our Cartesian coordinate system  is chosen in such a way that $x_0$  coincides with its origin (i.e. $x_0=0$) and the set $\Om(x_0, R)$  is described   by relations
\begin{equation}
\Om(x_0, R)\ = \ \Big\{  ~x= (x_1, x_2, x_3)\in \mathbb R^3~|~     x'\in S_R, \ \ph(x')< x_3<\ph(x')+\sqrt{R^2-|x'|^2} \ \Big\}.
\label{Om_R}
\end{equation}
Here we denote $x':=(x_1, x_2)$ and $S_R :=\{ ~x'\in \Bbb R^2~|~\sqrt{x_1^2 + x_2^2}<R~\}$.
With this assumptions the boundary condition in \eqref{Stokes_system-1} is equivalent to the relation
$$
u|_{x_3=\ph(x')} =0.
$$
We  assume $\ph$ is of class $W^3_\infty$ (i.e. its second derivatives are Lipschitz continuous) and the Cartesian coordinate system is chosen in such a way that the following relations hold
\begin{equation}
\ph(0)=0, \qquad \nabla \ph(0)=0, \qquad \| \ph\|_{W^3_\infty( S_R)}\le \mu.
\label{mu}
\end{equation}
Now we apply the diffeomorphism flattering the boundary, or, in other words, we introduce  new coordinates $y=\psi(x)$  by formulas
\begin{equation}
\psi: \Om_R \to B^+_R, \qquad
y= \psi(x) \ = \  \left( \begin{array}c x' \\ x_3-\ph(x') \end{array} \right),
\label{Definition_of_psi}
\end{equation}
$$
x\in \Om(x_0, R) \quad \Longleftrightarrow \quad y\in B_R^+.
$$
Denote
$$
\gathered
v:= u\circ \psi^{-1} , \qquad
q:= p\circ \psi^{-1}, \qquad \tilde f:=  f\circ \psi^{-1}.
\endgathered
$$
Then for $y=\psi(x)$ we have relations
$$
\gathered
\nabla p(x) = \hat \nabla_\ph q(y), \quad
\Delta u (x)  = \hat \Delta_\ph v (y), \quad
\div u(x) = (\hat\nabla_\ph \cdot v)(y).
\endgathered
$$
where $\hat \Delta_\ph$ and $\hat \nabla_\ph$ are  the differential operators
 with variable coefficients defined via a function $\ph$   by formulas
$$
\gathered
 \hat \Delta_\ph v \ := \ \Delta v- 2v_{, \al 3}\ph_{,\al} + v_{, 33} |\nabla'\ph|^2 - v_{,3} \Delta'\ph, \\
 \hat \nabla_\ph \cdot v \ := \ \div v - v_{\al, 3} \ph_{,\al}, \\
 \hat\nabla_\ph q \ := \  \nabla q - q_{,3} \left( \begin{array}c \nabla'\ph \\ 0 \end{array}\right).
 \endgathered
$$
  Here we assume summation from $1$ to $2$ over repeated Greek indexes and $\nabla'$ and $\Delta'$ denote the gradient and Laplacian with respect to $(y_1, y_2)$ variables.

The Navier-Stokes system \eqref{Stokes_system-1}  in $\Om(x_0, R)\times (-R^2,0)$ in $x$-variables is transformed to the  following system with variable coefficients depending on the  $y$-variables:
\begin{equation}
\left\{ \quad
\gathered
\gathered
\partial_t v   \ - \  \hat \Delta_\ph v \ + \ (v\cdot \hat \nabla_\ph ) v + \   \hat \nabla_\ph q \ = \ \tilde f
\\
\hat \nabla_\ph \cdot  v \ = \ 0
\endgathered
\qquad\mbox{in} \quad Q^+(R),
\\
v|_{y_3=0} = 0. \qquad\qquad\qquad\quad
\endgathered \right.
\label{Perturbed_Stokes}
\end{equation}
We  call this system  {\it the Perturbed Navier-Stokes system}. Note that if the boundary $\cd \Om$ is smooth in the neighborhood of $x_0$ then the coefficient of system \eqref{Perturbed_Stokes} are also smooth.

The Perturbed Navier-Stokes system possesses
the following scaling property:
if functions $v$, $q$, $f$, $\ph$ satisfy \eqref{Perturbed_Stokes} in the cylinder $Q^+(R)$ with $\ph$ satisfying \eqref{mu} then the functions
\begin{equation}
\gathered
v^R(x,t)= Rv(Rx, R^2t), \qquad q^R(x,t)=R^2q(Rx,R^2t), \\ f^R(x,t) = R^3f(Rx,R^2t), \qquad \ph^R(x')=\frac 1R\ph(Rx')
\endgathered
\label{Scaling}
\end{equation}
satisfy the Perturbed Navier-Stokes system \eqref{Perturbed_Stokes} in $Q^+$ and from Taylor decomposition of the function $\ph^R$ one can obtain for $R\le 1$
\begin{equation}
\ph^R(0)=0, \qquad \nabla'\ph^R(0) = 0, \qquad
\| \ph^R\|_{W^3_\infty( S_1)} \ \le \ \mu R.
\label{scaled_phi}
\end{equation}
Hence, if assumptions \eqref{mu} hold with some constant $\mu$ whose value can be arbitrary large (in particular, this implies that  the curvature of the boundary of $\Om$ in the neighborhood of $x_0\in \cd \Om$ can be arbitrary) then applying diffeomorphism \eqref{Definition_of_psi} and  the scaling transformation \eqref{Scaling} we can reduce the study of regularity of weak solutions to the Navier-Stokes system \eqref{Stokes_system-1} in domain $\Om(x_0, R)\times ]-R^2, 0[$  to the study of the Perturbed Navier-Stokes system \eqref{Perturbed_Stokes} in the canonical domain $Q^+$. Moreover, for any given $\mu_*>0$ choosing the radius $R=R(\mu_*)>0$  sufficiently small thanks to \eqref{scaled_phi} we can assume that condition
\begin{equation}
\mu R \le \mu_*
\label{mu_star}
\end{equation}
holds with some absolute  constant  $\mu_*>0$. If the value $\mu_*$ in \eqref{mu_star} is chosen sufficiently small then variable coefficients in the Perturbed Navier-Stokes system \eqref{Perturbed_Stokes} for functions $v^R$, $q^R$, $f^R$, $\ph^R$ in $Q^+$ can be interpreted as  small perturbations of the ``constant coefficients'' in the usual Navier-Stokes system \eqref{NSE_boubdary} for functions $v^R$, $q^R$, $f^R$ in $Q^+$.

The linear theory of strong solutions to the Stokes system extending  results of \cite{Seregin_ZNS271} to the case
of  the curvilinear  boundary  was developed in \cite{Solonnikov_ZNS288}. Later, in \cite{Shilkin_Vyalov} the similar theory was developed for the linear Perturbed Stokes system which is the linearization of \eqref{Perturbed_Stokes}. Moreover, in contrast to \cite{Solonnikov_ZNS288}, in \cite{Shilkin_Vyalov} the local estimates were obtained not for strong but for weak solutions. In particular, the analogs of our Theorems \ref{Theorem_Local_Estimate} --- \ref{Theorem2} were proved in \cite{Shilkin_Vyalov} for the Perturbed Stokes system under the assumption that its coefficients are  small perturbations of the constant coefficients of the usual Stokes system.

The linear theory developed in \cite{Solonnikov_ZNS288}, \cite{Shilkin_Vyalov} allows to prove the analogue of the basic $\ep$-regularity condition \eqref{Fixed_r_condition} at the neighborhood of point $x_0$ belonging to a curvilinear part of the boundary. To formulate a result we need to define   boundary suitable weak solution to the Navier-Stokes equation near a curvilinear part of the boundary. The definition of this class of solutions is analogues to one given in our Definition \ref{Definition_2} for the case of flat boundaries, see details in \cite{SSS}. So, we obtain the following result:

\begin{theorem}
\label{Main Theorem_Curved} Let $\Om(x_0, R)$ be defined by \eqref{Om_R} where $\ph\in W^3_\infty(S_R)$ satisfies \eqref{mu} and  assume $z_0=(x_0, t_0)$.  There exist an absolute constant  $\ep_*>0$ and a constant $R_*\in (0,R)$ depending only on $\mu$ and $R$
such that for any boundary suitable weak solution $u$ and $p$ of the
Navier-Stokes equations \eqref{Stokes_system-1} in $\Om(x_0, R)\times ]t_0-R^2, t_0[$   the following is true: if  there exists  $r \le R_*$ such that
\begin{equation}
\frac {1}{r^2} \
\int\limits_{t_0-r^2}^{t_0}\int\limits_{\Om(x_0,r)}
\left( |u|^3 + |p|^{\frac 32}\right)~dxdt < \ep_* \label{Rough}
\end{equation}
then $u$ is H\" older continuous on $ \bar\Om(x_0,\frac {r}2)
\times [t_0-\frac {r^2}{4}, t_0]$.
\end{theorem}

Theorem \ref{Main Theorem_Curved} is proved in \cite{SSS}. Using this  theorem one can obtain the following result that is a boundary version of the Caffarelli-Kohn-Nirenberg theorem, see \cite{CKN}:

\begin{theorem}
\label{CKN Theorem} Let all assumptions of Theorem \ref{Main Theorem_Curved} be satisfied. There exists an absolute constant $\ep_{**}>0$
such that, for any boundary suitable weak solution $u$ and $p$ of the
Navier-Stokes equations \eqref{Stokes_system-1} in $\Om(x_0, R)\times ]t_0-R^2, t_0[$, the velocity field $u$ is H\" older continuous in some neighborhood of $z_0$ provided
\begin{equation}
\limsup\limits_{r\to 0}~ \frac 1{r}
\int\limits_{t_0-r^2}^{t_0} \int\limits_{\Om(x_0,r)} |\nabla
u|^2~dxdt \ < \ \ep_{**}. \label{CKN_Condition_curv}
\end{equation}

\end{theorem}

Theorem \ref{CKN Theorem} provides the following estimate of the parabolic Hausdorff measure of singular set:

\begin{theorem}
\label{Estimate Sigma} Let $\Om\subset \Bbb R^3 $ be a domain whose boundary  $\cd \Om$ is of class $W^3_\infty$ and assume there is a constant $\mu>0$ such that for any $x_0\in \cd \Om$ there exists a neighborhood $\Om(x_0, R)$ which can be described in an appropriate Cartesian coordinate system by formulas \eqref{Om_R} with some function $\ph=\ph_{x_0}$ satisfying conditions \eqref{mu}.
Then for any boundary suitable weak solution $u$ and $p$ of the Navier-Stokes system in $\Om \times ]0,T[$ there exists a closed set
$\Sigma \subset \cd \Om\times ]0,T]$ such that for any point
$z_0\in (\cd \Om\times ]0,T]) \setminus \Sigma$ the function $u$
is H\" older continuous in some neighborhood of $z_0$ and,
moreover,
$$
\mathcal P^1(\Sigma) =0,
$$
where $\mathcal P^1(\Sigma)$ is the one-dimensional parabolic Hausdorff measure of $\Sigma$.

\end{theorem}
\medskip

 It worthy to notice that there is an essential difference between $\ep$-regularity conditions \eqref{Rough} and \eqref{CKN_Condition_curv}. Condition \eqref{CKN_Condition_curv} requires smallness of functional
 $$
 E(u,r) = \Big( \frac 1r ~\int\limits_{t_0-r^2}^{t_0} \int\limits_{\Om(x_0,r)} |\nabla
u|^2~dydt\Big)^{1/2}
 $$
 uniformly with respect to all $r\in(0,R_1)$ with some $R_1>0$, i.e.
\begin{equation}
\limsup\limits_{r\to 0} E(u, r) < \ep_0.
\label{limsup_condition}
\end{equation}
 The question is whether (\ref{limsup_condition}) could be weakened to   
\begin{equation}
\liminf\limits_{r\to 0} E(u, r) < \ep_0.
\label{limsup_condition_1}
\end{equation}
 In general, an answer to this question is unknown. 
But  if we assume {\it a priori} boundedness of scale invariant functionals of type \eqref{Functionals} then
it is possible to improve most of $\ep$-regularity conditions, replacing in \eqref{limsup_condition} $\limsup$ with $\liminf$. Namely,  the boundary suitable weak solution $u$ and $p$ to the Navier-Stokes system in $Q^+$ is regular near the origin if \eqref{Boundedness} holds and besides one of the following conditions is valid:
\begin{itemize}
\item $\min\Big\{~\liminf\limits_{r\to 0} A(u,r), ~ \liminf\limits_{r\to 0} C(u,r), ~\liminf\limits_{r\to 0} E(u,r), ~\liminf\limits_{r\to 0} D(p,r)  ~\Big\} <\ep_0$
\item $\liminf\limits_{r\to 0} C_2(u,r)  <\ep_0$, \ where
$$
C_2(u,r) := \Big(~\frac 1{r^3}~\int\limits_{Q^+(r)} |u(x,t)|^2~dxdt~\Big)^{1/2}
$$
\end{itemize}
This statement has been  proven  in  \cite{Seregin_UMN} in the internal case. Later on, it was extended to the boundary case in  \cite{Mikhaylov_1}, \cite{Mikhaylov_2}.

\end{document}